\documentclass[review]{elsarticle}

\usepackage{amssymb,amsmath}
\usepackage[usenames,dvipsnames]{xcolor}
\usepackage{subfigure}
\usepackage{caption}
\usepackage{algorithm, algorithmicx, algcompatible}
\usepackage{graphicx}
\usepackage{tikz}
\usepackage[export]{adjustbox}
\usepackage{epstopdf}

\definecolor{cadmiumgreen}{rgb}{0.0, 0.42, 0.24}

\usepackage{todonotes}

\usepackage{lineno}
\modulolinenumbers[5]

\newdefinition{rmk}{Remark }

\makeatletter
\def\@author#1{\g@addto@macro\elsauthors{\normalsize%
    \def\baselinestretch{1}%
    \upshape\authorsep#1\unskip\textsuperscript{%
      \ifx\@fnmark\@empty\else\unskip\sep\@fnmark\let\sep=,\fi
      \ifx\@corref\@empty\else\unskip\sep\@corref\let\sep=,\fi
      }%
    \def\authorsep{\unskip,\space}%
    \global\let\@fnmark\@empty
    \global\let\@corref\@empty  
    \global\let\sep\@empty}%
    \@eadauthor={#1}
}
\makeatother

\begin{document}

\begin{frontmatter}

\title{The Multiscale Perturbation Method for Two-Phase Reservoir Flow Problems}

\author[1]{Franciane F. Rocha}
\cortext[cor1]{Corresponding author}
\ead{fr.franciane@usp.br}
\author[2]{Het Mankad}
\ead{Het.Mankad@utdallas.edu}
\author[1]{Fabricio S. Sousa}
\ead{fsimeoni@icmc.usp.br}
\author[2]{Felipe Pereira\corref{cor1}}
\ead{luisfelipe.pereira@utdallas.edu}

\address[1]{Instituto de Ci\^encias Matem\'aticas e de Computa\c c\~ao, \\
Universidade de S\~ao Paulo, Av. Trabalhador S\~ao-carlense, 400,  \\
13566-590, S\~ao Carlos, SP, Brazil}

\address[2]{Department of Mathematical Sciences, The University of Texas at Dallas, \\
800 W. Campbell Road, Richardson, Texas 75080-3021, USA}


\begin{abstract}
In this work we formulate and test a new procedure, the Multiscale Perturbation Method for Two-Phase Flows (MPM-2P), for the fast, accurate and naturally 
parallelizable numerical solution of two-phase, incompressible, immiscible displacement in porous media approximated by an operator splitting method.
The proposed procedure is based on domain decomposition and combines the Multiscale 
Perturbation Method (MPM) [Ali, et al., \textit{Appl. Math. and Comput.},  125023 (2020)] with the Multiscale Robin Coupled Method (MRCM) 
[Guiraldello, et al., \textit{J. Comput. Phys.},  355 (2018) pp. 1-21]. 
When an update of the velocity field is called by the operator splitting algorithm, 
 the MPM-2P may provide, depending on the magnitude of a dimensionless algorithmic parameter, an accurate and computationally inexpensive 
 approximation for the velocity field by reusing previously computed multiscale basis functions. Thus, a full update of all multiscale basis functions 
 required by the MRCM for the construction of a new velocity field is avoided. 

There are two main steps in the formulation of the MPM-2P.
Initially, for each subdomain one local boundary value problem with trivial Robin boundary conditions is solved (instead of 
a full set of multiscale basis functions, that would be required by the MRCM). Then, the solution of an inexpensive interface problem provides
the velocity field on the skeleton of the decomposition of the domain. The resulting approximation for the velocity field is obtained by downscaling.

We consider challenging two-phase flow problems, with high-contrast permeability fields and water-oil finger growth in homogeneous media. 
Our numerical experiments show that the use of the MPM-2P gives exceptional  speed-up - almost 90\% of reduction in computational cost - of two-phase flow simulations. 
Hundreds of MRCM solutions can be replaced by inexpensive MPM-2P solutions, and water breakthrough can be simulated with very few updates of the MRCM set of multiscale basis functions. 
\end{abstract}
\end{frontmatter}
\section{Introduction}

We are concerned with the development of fast and scalable multiscale solvers for porous media flows, aiming at the solution of inverse problems (uncertainty
quantification) in reservoir simulation, where thousands of forward in time simulations have to be performed (see, for instance \cite{ALMAMUN2020124980} and references
therein). In this work we consider two-phase, incompressible, immiscible displacement in porous media approximated by an operator splitting method. Such procedure
decomposes the governing coupled system of partial differential equations into two equations that are solved sequentially: a second 
order elliptic equation for the fluid pressure (the pressure equation) is followed by the solution of a scalar hyperbolic conservation law for a phase saturation
(the saturation equation). Our focus here is to speed-up the solution of the pressure equation by combining recent developments in the area of multiscale mixed methods for second order elliptic equations. 
 Within the splitting framework the hyperbolic conservation law is solved by an explicit finite volume scheme, that
can be efficiently solved in multi-core devices.

There are many different types of multiscale methods that can be used to solve the elliptic equation frequently occurring in problems related to the reservoir simulations.  These methods are typically based on the finite element method (FEM), finite volume method (FVM) or mixed finite element method. A detailed review of some of the established multiscale methods based on these various procedures can be found in \cite{Kippe}. 
Method discussed in \cite{HouMulti,aarnes2002multiscale} and \cite{efendiev2006accurate} are FEM based while \cite{jenny2003multi} and \cite{jenny2005adaptive} are FVM based multiscale methods used to solve the second order elliptic equation. There are also other different approaches like the variational multiscale method \cite{hughes1998variational}, the generalized multiscale method \cite{efendiev2013generalized} and the hybrid multiscale method \cite{madureira2017hybrid}. We are mainly concerned about the multiscale methods based on mixed FEM approach (see \cite{ARNOLD1990281} for  details about the theoretical aspect of the mixed FEM). These type of methods can be found, for example, in \cite{aarnes,Arbogast::2007,pereira,chen_hou}. Recent developments in this area of research have been introduced in \cite{guiraldello2018multiscale} and \cite{guiraldello2019interface}.    

We have recently established in \cite{Paolathesis,PaolaJCP} that a recursive formulation of the Multiscale Robin Coupled Method 
(MRCM) \cite{guiraldello2018multiscale} shows excellent scalability (both weak and strong) for the solution of the pressure equation. 
These conclusions were reached by solving the pressure equation on state-of-the-art multi-core devices, for problems with a few billion variables,
 that are of interest to the oil industry. 
In \cite{Paolathesis,PaolaJCP} the solution of a second order elliptic equation is obtained in two steps.
In a first step, for each subdomain of a decomposition of the domain of interest a set of multiscale basis functions (local boundary value problems of Robin type)
has to be computed. Then, a coarse interface problem defined on the skeleton of the domain decomposition needs to be solved. It has been shown 
in \cite{Paolathesis,PaolaJCP} that the time associated with the solution of the interface problem is essentially {\it negligible}, when compared to the time spent in solving 
the local boundary value problems that give the multiscale basis functions. Thus, a fair assessment of the cost of the solution of the pressure equation
by a multiscale method can be made in terms of the number of multiscale basis functions that are computed. Our main objective in this work is to design a 
method that can accomplish a reduction in the number of updates of multiscale basis functions that are needed in the numerical simulation of two-phase flows.

In order to explain our approach, consider the operator splitting scheme for two-phase flows \cite{douglas1983time, douglas1997numerical, furtado2011operator}. If the above mentioned multiscale mixed method is applied to solve 
the pressure equation, then a set of multiscale basis functions has to be, in principle, recomputed every time the solution algorithm calls for an updated velocity 
field. Thus, the development of methods that reduce the number of multiscale basis functions to be computed in each subdomain, without loss of accuracy, is
of great importance to speed-up the solution of two-phase flow problems. The procedure that we introduce in this work, the Two-Phase Multiscale Perturbation Method
(MPM-2P) has precisely this objective. The MPM-2P is based on the original Multiscale Perturbation Method (MPM) ~\cite{ALI2020125023}, that was introduced 
to approximate the velocity field by reusing multiscale basis functions computed for a distinct pressure equation (with different, but closely related coefficients), 
provided that the solutions of the two elliptic equations at hand can be related by classical perturbation theory \cite{nayfeh2011introduction}. The proposed method combines the MPM
with the most recent developments of the MRCM for two-phase flow problems \cite{bifasico}. Our results are very encouraging. We consider challenging two-phase flow simulations and we find that we can typically reduce the computational 
cost of a simulation by up to two orders of magnitude. 
Although the new method is presented for two-phase flows it can also be applied to other types of flows,
as well as to the sequential implicit solution of multiphase flows. These developments, including the implementation of the procedure described here in multi-core and 
multi-GPU devices, are currently being considered by the authors and their collaborators.

This work is organized as follows. The model equations for two-phase flows appear in Section 2. In Section 3, we recall the MRCM and present a description of special multiscale basis functions that will be used in our experiments. Then, the new algorithm for the operator splitting method based on perturbation techniques  which we call the MPM-2P is introduced in Section 4. In this section, at first we will briefly discuss the MPM-2P and then we propose a modified operator splitting scheme. A cost analysis of the new method is presented which shows its computational efficiency when compared to a classic operator splitting scheme. Our numerical experiments are presented in Section 5, followed by Section 6 with the concluding remarks.


\section{The Two-Phase Flow Problem} 
In this section we will introduce the model equations for the two-phase flow problem where the two phases into consideration are oil and water
 (see ~\cite{zchen_book,ewing1983mathematics} for a discussion of these equations). The governing system of equations that we will describe here is in a dimensionless form (see  ~\cite{ALI2020125023} for a detailed description of the dimensionless form). The system is given by the following elliptic problem 
\begin{equation}\label{Darcy} 
\begin{array}{rll}
\mathbf{u}&=-\lambda(s)K(\mathbf{x})\nabla p  &\mbox{in}\ \Omega ,\\
\nabla \cdot \mathbf{u}&=q  &\mbox{in}\ \Omega, \\ 
p &= g &\mbox{on}\ \partial\Omega_{p},\\
\mathbf{u} \cdot \mathbf{n}&= z &\mbox{on}\ \partial\Omega_{u},
\end{array}
\end{equation}
coupled with the hyperbolic conservation law for the saturation transport problem
\begin{equation}\label{BL2D}
\begin{array}{rll}
\dfrac{\partial s}{\partial t} + \nabla \cdot \left(f(s)\mathbf{u}\right)& = 0  &\mbox{in}\ \Omega, \\
s(\mathbf{x},t=0) &= s^0(\mathbf{x}) &\mbox{in}\ \Omega,\\
s(\mathbf{x},t) &= \bar{s}(\mathbf{x},t) &\mbox{in}\ \partial\bar{\Omega}.
\end{array}
\end{equation}
Here the velocity $\mathbf u (\mathbf x, t)$, pressure $p(\mathbf x,t)$ fields and the water saturation 
$s(\mathbf x, t)$ are the unknown quantities of the model problem in the domain $\Omega$. We consider the 2D system for our current experiments, but all developments here can be extended to 3D space as well, without any additional effort. In this work, capillary pressure and gravity effects are not taken into account, but they can be added without changing the proposed method. Moreover, we also consider a fully saturated media which is a common practice in problems related to oil-water flows in petroleum reservoirs, which means that $s_w + s_o = 1$, hence the saturation equation needs to consider only one of them, say $s := s_w$ as stated in Eq. (\ref{BL2D}). Here, $K(\mathbf{x})$ is the absolute permeability; $q=q(\mathbf x, t)$ is a source term; $g = g(\mathbf x, t)$ is the boundary condition for pressure at $\partial\Omega_{p}$; $z = z(\mathbf x, t)$ is the boundary condition for normal component of flux at $\partial\Omega_{u}$ ($\mathbf{n}$ is the outward unit normal); $s^0(\mathbf{x})$ is the saturation initial condition; and $\bar{s}(\mathbf{x},t)$ is the saturation at the injection boundaries $\partial\Omega^-=\{\mathbf{x}\in\partial\Omega,\ \mathbf{u} \cdot \mathbf{n}<0\}$. The coefficient $\lambda(s) = \lambda_{o}(s) + \lambda_{w}(s)$ is the total phase mobility where,
\begin{align}
\lambda_{j}(s) = \frac{k_{rj}(s)}{\mu_{j}},
\end{align}
with  $j \in \{w,\ o\}$ representing the water and the oil phase, while ${k}_{rj}(s)$ and $\mu_{j}$ are the relative permeability field and viscosity corresponding to the $j$ phase respectively.  The fractional flow function of water is given as,
\begin{align}
f(s) = \frac{\lambda_w(s)}{\lambda(s)}.
\end{align}
Our model considers a constant porosity scaled out by changing the time variable. Henceforth, we will denote the conductivity by $\kappa(\mathbf x) = \lambda(s) K(\mathbf x)$. 

\section{The Multiscale Robin Coupled Method} 

In this section we will provide a brief overview of the MRCM \cite{guiraldello2018multiscale} which is also a generalization of  the Multiscale Mixed Method (MuMM) \cite{pereira}.  MRCM is a non-overlapping domain decomposition method in that the domain $\Omega$ is divided into several subdomains $\Omega_i, i=1, 2, \cdots, N$. The MRCM algorithm to approximate the solution of the second order elliptic equation consists of two parts. In the first part, the  solution to the elliptic equation is approximated locally for each of the subdomains $\Omega_i,  i=1,2,\cdots, N$. In the second part, the coarse interface problem defined on the skeleton $\Gamma$ of the domain decomposition (the union of all interfaces $\Gamma_{ij} = \Omega_i \cap \Omega_j$) has to be solved.
We refer to two different length scales: $h$, the fine mesh size, and $H$, the characteristic size of the subdomains also known as coarse scale. Here, $H\gg h$. We will denote the elliptic solution  obtained by MRCM for each fine scale element by $(\mathbf{u}_h, p_h)$.

Weak continuity of the solution is imposed on the coarse scale through the following compatibility conditions:
\begin{equation}
\int_\Gamma (\mathbf u_h^+ - \mathbf u_h^-)\cdot \check{\mathbf n} \ \psi \ d\Gamma=0 \quad \text{and}\quad 
\int_\Gamma (p_h^+ - p_h^-)\ \phi \ d\Gamma=0.
\label{eq:compatibility}
\end{equation}
Here $(\phi,\psi) \in \mathcal{U}_H \times \mathcal{P}_H$ where $\mathcal{U}_H \text{ and } \mathcal{P}_H$ are low-dimensional interface spaces defined over the edges $\mathcal{E}_h$ of the skeleton $\Gamma$, that are subspaces of
\begin{equation}
\mathfrak{F}_h(\mathcal{E}_h) = \left\{ f:\mathcal{E}_h\to \mathbb{R};~f|_e\,\in\,\mathbb{P}_0~,~\forall\,e\,\in\,\mathcal{E}_h \right\}.
\end{equation}
In Eq. (\ref{eq:compatibility}), the solution on each side of the interface $\Gamma$ is represented by the $+$ and $-$ superscripts, while the normal vector to the skeleton is denoted by $\check{\mathbf n}$. 
These compatibility conditions are enforced by imposing the following Robin-type boundary conditions to the local problems
\begin{equation}
-\beta_i  \mathbf u_h^i \cdot \check{\mathbf n}^i + p_h^i =  -\beta_i U_H \check{\mathbf n} \cdot \check{\mathbf n}^i + P_H 
\end{equation}
where $(\mathbf u_h^i, p_h^i)$ are the local normal flux and pressure unknowns for each subdomain $\Omega_i$ and $(U_H,P_H)$ are global unknowns defined on the interface of the decomposition of the domain.
The parameter $\beta_i$ on each subdomain is defined as
\begin{equation}
\beta_i(\mathbf {x}) = \frac{\alpha(\mathbf{x}) H}{\kappa_i(\mathbf x)},
\label{eq:beta}
\end{equation}
where $\alpha(\mathbf {x})$ is a dimensionless algorithmic function that is locally defined according to the variations in the permeability field (see \cite{bifasico}).

The MRCM is formulated as : Find $(\mathbf u_h^i, p_h^i)$ and $(U_H,P_H)$ such that the following local problems are satisfied
\begin{equation}
\begin{array}{rcll}
\mathbf u_h^i &=& -\kappa (\mathbf x) \ \nabla p_h^i & \text{in }\Omega_i, \\
\nabla \cdot \mathbf u_h^i &=& q & \text{in }\Omega_i, \\
p_h^i & = & g_p & \text{on }\partial\Omega_i\cap \partial\Omega_p, \\
\mathbf u_h^i \cdot \check{\mathbf n}^i &=& g_u & \text{on }\partial\Omega_i\cap \partial\Omega_u, \\
-\beta_i  \mathbf u_h^i \cdot \check{\mathbf n}^i + p_h^i &=&  -\beta_i U_H \check{\mathbf n} \cdot \check{\mathbf n}^i + P_H & \text{on }\partial\Omega_i\cap \Gamma, 
\end{array}
\label{eq:mrc1}
\end{equation}
along with the following global system
\begin{equation}
\begin{array}{rcl}
\displaystyle \sum_{i=1}^N \int_{\partial\Omega_i\cap \Gamma} (\mathbf u_h^i \cdot \check{\mathbf n}^i) \ \psi \ d\Gamma &=& 0,\\
\displaystyle \sum_{i=1}^N \int_{\partial\Omega_i\cap \Gamma} \beta_i (\mathbf u_h^i \cdot \check{\mathbf n}^i - U_H\ \check{\mathbf n} \cdot \check{\mathbf n}^i)\ \phi  \ (\check{\mathbf n} \cdot \check{\mathbf n}^i)\ d\Gamma &=& 0,
\end{array}
\label{eq:mrc2}
\end{equation}
for all $(\phi,\psi)\in \mathcal{U}_H\times\mathcal{P}_H$.

The implementation of the MRCM considers an additive decomposition of the local solutions $(\mathbf u_h^i, p_h^i)$ given by
\begin{equation}
\mathbf{u}_h^i=\hat{\mathbf{u}}_h^i+\bar{\mathbf{u}}_h^i, \qquad p_h^i=\hat{p}_h^i+\bar{p}_h^i,
\end{equation}
that satisfies
\begin{equation}
\begin{array}{rcll}
\hat{\mathbf u}_h^i &=& -\kappa (\mathbf x) \ \nabla \hat{p}_h^i & \text{in }\Omega_i \\
\nabla \cdot \hat{\mathbf u}_h^i &=& 0 & \text{in }\Omega_i \\
\hat{p}_h^i & = & 0 & \text{on }\partial\Omega_i\cap \partial\Omega_p \\
\hat{\mathbf u}_h^i \cdot \check{\mathbf n}^i &=& 0 & \text{on }\partial\Omega_i\cap \partial\Omega_u \\
-\beta_i  \hat{\mathbf u}_h^i \cdot \check{\mathbf n}^i + \hat{p}_h^i &=&   -\beta_i U_H \check{\mathbf n} \cdot \check{\mathbf n}^i + P_H & \text{on }\partial\Omega_i\cap \Gamma 
\end{array}
\label{eq:mrc_homogeneous}
\end{equation}
and
\begin{equation}
\begin{array}{rcll}
\bar{\mathbf u}_h^i &=& -\kappa (\mathbf x) \ \nabla \bar{p}_h^i & \text{in }\Omega_i \\
\nabla \cdot \bar{\mathbf u}_h^i &=& q & \text{in }\Omega_i \\
\bar{p}_h^i & = & g_p & \text{on }\partial\Omega_i\cap \partial\Omega_p \\
\bar{\mathbf u}_h^i \cdot \check{\mathbf n}^i &=& g_u & \text{on }\partial\Omega_i\cap \partial\Omega_u \\
-\beta_i  \bar{\mathbf u}_h^i \cdot \check{\mathbf n}^i + \bar{p}_h^i &=& 0 & \text{on }\partial\Omega_i\cap \Gamma .
\end{array}
\label{eq:mrc_non_homogeneous}
\end{equation}
The local problems in Eq. (\ref{eq:mrc_homogeneous}) satisfy a nonzero Robin boundary condition for the subdomain coupling and have source terms as well as physical boundary conditions identically equal to zero. They correspond to the homogeneous part of the solution. 
On the other hand, the local problems in Eq. (\ref{eq:mrc_non_homogeneous}) have interface Robin boundary condition for the subdomain coupling equal to zero and take into account the contribution of the nonzero source terms as well as the nonzero boundary conditions. This solution represents the non-homogeneous part of the final numerical solution. The set of solutions generated by numerically solving Eq. (\ref{eq:mrc_homogeneous}) forms a set of multiscale basis functions (BFs). Moreover,  the solution to Eq. (\ref{eq:mrc_non_homogeneous}) gives one additional local BF.

The interface spaces $\mathcal{U}_{H} \text{ and } \mathcal{P}_{H}$ are spanned by the multiscale BFs \\ $\{\phi_1, \phi_2, \cdots , \phi_{N_U}\}$ and $\{\psi_1, \psi_2, \cdots ,   \psi_{N_P}\}$, where $N_U = \dim (\mathcal{U}_H)$ and $N_P = \dim (\mathcal{P}_H)$. Thus, the interface unknowns $U_H$ and $P_H$ are given by,
\begin{equation}
U_H=\displaystyle \sum_{l=1}^{N_U} U_l\phi_l,\qquad P_H=\displaystyle \sum_{l=1}^{N_P} P_l\psi_l,
\end{equation}
where the coefficients $U_l$ and $P_l$ are the solution of the global interface system generated by the Eq. (\ref{eq:mrc2}) when tested with all BFs of $\mathcal{U}_H$ and $\mathcal{P}_H$. We remark that the local problems are completely independent and can be computed in parallel. 

\subsection{Choice of interface spaces for the multiscale basis functions}\label{sec:mrcm_bf}

 In this section we will discuss our choice of the interface spaces $\mathcal{U}_H \text{ and } \mathcal{P}_H$ that we will be using to obtain the solution of the multiscale BFs. In terms of degrees of freedom per interface,  $N_U=k_U\times N_I$ and $N_P = k_P\times N_I$, where $k_U$, $k_P$ and $N_I$ are, respectively, the flux degrees of freedom, pressure degrees of freedom and number of interfaces between subdomains (see \cite{guiraldello2018multiscale} for more details). 

The use of classic low-degree polynomial functions (projected onto $\mathfrak{F}_h(\mathcal{E}_h) $) is the most common choice for the interface spaces $\mathcal{U}_H$ and $\mathcal{P}_H$. For Gaussian permeability fields, these spaces are enough to ensure accurate approximations by choosing linear polynomial interfaces. However, for high-contrast channelized permeability fields, such as the ones considered here, polynomial based spaces are not adequate to capture these types of features. Alternatives are informed spaces, as in 
\cite{guiraldello2019interface}, or the use of recently developed spaces based on physics \cite{rocha2020interface, rochaenhanced}, which are capable of accurately capturing homogeneities such as channels and barriers, as happens in fractured karstified reservoirs \cite{popov2009multiphysics, lopes2019new}.

In short the new interface space based on physics is an adaptive piecewise polynomial (further projected onto $\mathfrak{F}_h(\mathcal{E}_h) $) that automatically accomodates pressure discontinuities across high permeability channels, as well as flux discontinuities across low-permeability barriers. Such adaptive spaces are capable of recovering the true physical solution of the flow in presence of these heterogeneities.

This strategy is seamlessly combined with the adaptive version of the MRCM (called $a$MRCM) to set values of $\alpha(\mathbf{x})$ function according to permeability variations \cite{bifasico}. According to the authors, the $a$MRCM is able to reduce the error introduced by the domain decomposition if small values of $\alpha$ are chosen for high permeable regions, whilst large values are chosen for remaining areas. The combination of $a$MRCM \cite{bifasico} with interface spaces based on physics, as shown in \cite{guiraldello2019interface}, seems to be the most accurate strategy to deal with highly heterogeneous media and therefore this is our choice of multiscale domain decomposition method for the numerical simulations presented in this work.

\section{The Multiscale Perturbation Method for Two-Phase Flows}


We consider an operator splitting scheme for two-phase flows as presented in \cite{bifasico}, where pressure and saturation are updated sequentially (see \cite{douglas1983time, douglas1997numerical, furtado2011operator} for additional discussions about the operator splitting framework).  
The pressure is updated at times $t_n = n\Delta t_p$, for $n=0,1,\dots$, while the saturation is computed at intermediate times $t_{n,k} = t_n + k\Delta t_s$, for $k=1,2,\dots, C_n$, such that $t_n < t_{n,k} \leq t_{n+1}$. Here, $\Delta t_s$ denotes the time step used in the discretization of the saturation equation, $\Delta t_p$ is the time step for pressure, and $C_n$ is the number of transport time steps between $t_n$ and $t_{n+1}$.

Let $p^n (\mathbf x)$, $\mathbf u^n (\mathbf x)$ and $s^n (\mathbf x)$ denote the pressure, velocity and saturation approximations at time $t^n$. We compute the saturation $s^{n} (\mathbf x)$ through Eq. (\ref{BL2D}) by using an explicit Euler time integration (with $\mathbf{u}^{n-1}$ constant at intermediate times $t_{n-1,k}$) combined with a first order upwind method \cite{leveque2002finite}. Then, the saturation $s^{n} (\mathbf x)$ is used to compute the pressure $p^{n} (\mathbf x)$ and velocity $\mathbf u^{n} (\mathbf x)$ through Eq. (\ref{Darcy}) by applying a multiscale method. At this point, instead of calling directly the MRCM, our operator splitting algorithm uses it in the framework of the MPM.

\subsection{Reusing previously computed basis functions}
The goal of the MPM-2P is to approximate the pressure $p^{n} (\mathbf{x})$ and velocity $\mathbf{u}^{n} (\mathbf{x})$ by reusing the BFs that are computed by the MRCM at an earlier time of the simulation. In order to introduce the formulation of the MPM-2P, consider that the BFs computed at time $t_m$ ($m<n$) will be reused.
Therefore, we have two elliptic problems: $\mathcal{P}_{t_m}$ and $\mathcal{P}_{t_n}$, associated with times $t_m$ and $t_n$, respectively. Following the perturbation theory presented in \cite{ALI2020125023}, we express the conductivity for problem $\mathcal{P}_{t_n}$ as a perturbation of the conductivity of problem $\mathcal{P}_{t_m}$, i.e. $\kappa^n=\kappa^n(\mathbf x)  = \lambda(s^n(\mathbf x)) K(\mathbf x)= \kappa^m+\epsilon\kappa_\epsilon$, where $\epsilon=||\kappa^n-\kappa^m||$ is a small parameter that measures the difference (in $L^2$ norm) between $\kappa^n$ and $\kappa^m$ from times $t_n$ and $t_m$, respectively, while $\kappa_\epsilon=(\kappa^n-\kappa^m)/\epsilon$ is an auxiliary field related to the formulation of the MPM. 
Thus, the two elliptic problems at hand can be written as
\begin{equation}\label{eq:mpm_t0}
\mathcal{P}_{t_m} :\ \left\{
\begin{array}{rll}
\mathbf{u}^m &= - \kappa^m\nabla p^m &\mbox{in}\ \Omega \\
\nabla \cdot \mathbf{u}^m &= q &\mbox{in}\ \Omega \\
p^m &= g &\mbox{on}\ \partial\Omega^{p} \\
\mathbf{u}^m \cdot \mathbf{n}&= z &\mbox{on}\ \partial\Omega_{u}
\end{array} \right. 
\end{equation}
and
\begin{equation}\label{eq:mpm_tf}
\mathcal{P}_{t_n} : \left\{
\begin{array}{rll}
\mathbf{u}^n &= - (\kappa^m+\epsilon\kappa_\epsilon)\nabla p^n &\mbox{in}\ \Omega \\
\nabla \cdot \mathbf{u}^n &= q &\mbox{in}\ \Omega \\
p^n &= g &\mbox{on}\ \partial\Omega_{p} \\
\mathbf{u}^n \cdot \mathbf{n}&= z &\mbox{on}\ \partial\Omega_{u}.
\end{array} \right. 
\end{equation}
For simplicity, we assume that the source term and the known boundary functions depend only on space, but time-dependent source terms and boundary data can, in principle, be considered.
Next, we write the pressure and flux of problem $\mathcal{P}_{t_n}$ as perturbations of the respective pressure and flux of problem $\mathcal{P}_{t_m}$:
\begin{equation}\label{eq:soma_p}
    p^n=p^m+ \delta p^n,     
\end{equation}
\begin{equation}\label{eq:soma_u}
     \mathbf{u}^n=\mathbf{u}^m+\delta\mathbf{u}^n.
\end{equation}
By combining this decomposition with problem $\mathcal{P}_{t_n}$ (\ref{eq:mpm_tf}) we get the following auxiliary system for the pair $(\delta\mathbf{u}^n,\delta p^n)$ 
\begin{equation}\label{eq:mpm_tf_2}
\mathcal{P}_{\hat{\mathbf{u}}} : \left\{
\begin{array}{rll}
\hat{\mathbf{u}} &= - (\kappa^m+\epsilon\kappa_\epsilon)\nabla \delta p^n &\mbox{in}\ \Omega \\
\nabla \cdot \hat{\mathbf{u}} &= q + \nabla \cdot((\kappa^m+\epsilon\kappa_\epsilon)\nabla p^m)&\mbox{in}\ \Omega \\
\delta p^n &= g -p^m &\mbox{on}\ \partial\Omega_{p} \\
\hat{\mathbf{u}} \cdot \mathbf{n}&= z + ((\kappa^m+\epsilon\kappa_\epsilon)\nabla p^m) \cdot \mathbf{n} &\mbox{on}\ \partial\Omega_{u},
\end{array} \right. 
\end{equation}
where $\hat{\mathbf{u}} = \mathbf{u}^m + \delta \mathbf{u}^n + (\kappa^m+\epsilon\kappa_\epsilon)\nabla  p^m$. Although this system is well-posed, its solution is as expensive as the direct solution of problem $\mathcal{P}_{t_n} $, so approximations are needed in order to reduce the cost of solving this auxiliary problem.

Since our goal is to reuse the BFs computed for problem $\mathcal{P}_{t_m}$ (\ref{eq:mpm_t0}), we
need to somehow connect the solution of $\mathcal{P}_{\hat{\mathbf{u}}}$ (\ref{eq:mpm_tf_2}) to the solution space of $\mathcal{P}_{t_m}$. This would allow us to write the solution of $\mathcal{P}_{\hat{\mathbf{u}}}$ by taking advantage of the span of the precomputed BFs. To approximate $\delta p^n$ and $\delta\mathbf{u}^n$, we consider the following perturbation expansions:
\begin{equation}
    \delta p^n=\delta p_0^n + \epsilon\delta p_1^n +\epsilon^2\delta p_ 2^n +  \epsilon^3\delta p_ 3^n +\cdots
    \label{eq:delta_p0}
\end{equation}
and
\begin{equation}
    \delta \mathbf{u}^n=\delta \mathbf{u}_0^n + \epsilon\delta \mathbf{u}_1^n +\epsilon^2\delta \mathbf{u}_ 2^n +  \epsilon^3\delta \mathbf{u}_ 3^n +\cdots 
     \label{eq:delta_u0}
\end{equation}
By applying Eqs. (\ref{eq:delta_p0}) and (\ref{eq:delta_u0}) in Eq. (\ref{eq:mpm_tf_2}), and considering the expansions up to term $\epsilon^{\ell}$, we get the following problems for $\ell=0$ and $\ell>0$, respectively:
\begin{equation}\label{eq:mpm_tf_3}
\mathcal{P}_{\hat{\mathbf{u}}_{0}} : \left\{
\begin{array}{rll}
\hat{\mathbf{u}}_{0} &= - \kappa^m\nabla \delta p_0^n &\mbox{in}\ \Omega \\
\nabla \cdot \hat{\mathbf{u}}_{0} &= q + \nabla \cdot((\kappa^m+\epsilon\kappa_\epsilon)\nabla p^m)&\mbox{in}\ \Omega \\
\delta p_0^n &= g -p^m &\mbox{on}\ \partial\Omega_{p} \\
\hat{\mathbf{u}}_{0} \cdot \mathbf{n}&= z + ((\kappa^m+\epsilon\kappa_\epsilon)\nabla p^m) \cdot \mathbf{n} &\mbox{on}\ \partial\Omega_{u},
\end{array} \right. 
\end{equation}
\begin{equation}\label{eq:mpm_tf_4}
\mathcal{P}_{\hat{\mathbf{u}}_{\ell}} : \left\{
\begin{array}{rll}
\hat{\mathbf{u}}_{\ell} &= - \kappa^m\nabla \delta p_{\ell}^n &\mbox{in}\ \Omega \\
\nabla \cdot \hat{\mathbf{u}}_{\ell} &= \nabla \cdot(\kappa_\epsilon\nabla p_{\ell-1}^n)&\mbox{in}\ \Omega \\
\delta p_{\ell}^n &= 0 &\mbox{on}\ \partial\Omega_{p} \\
\hat{\mathbf{u}}_{\ell} \cdot \mathbf{n}&= 0 &\mbox{on}\ \partial\Omega_{u},
\end{array} \right. 
\end{equation}
where $\hat{\mathbf{u}}_{0} = \mathbf{u}^m + \delta \mathbf{u}_0^n + (\kappa^m+\epsilon\kappa_\epsilon)\nabla  p^m$ and  $\hat{\mathbf{u}}_{\ell} = \delta\mathbf{u}_{\ell}^n +\kappa_\epsilon\nabla \delta p_{\ell-1}^n$.  
The problem $\mathcal{P}_{\hat{\mathbf{u}}_{0}}$ is associated with $\epsilon^0$, while the problems $\mathcal{P}_{\hat{\mathbf{u}}_{\ell}}$ are associated with $\epsilon^{\ell}$, for each $\ell>0$. 
\par
The approach developed in \cite{ALI2020125023}, neglects the perturbation terms with $\ell>1$, since they are small enough ($\epsilon^2\ll 1$ if the perturbation in the conductivity is small enough from time $t_m$ to time $t_n$). We find from our numerical experiments for two-phase flows that, even the first-order perturbation term can be neglected, still resulting in accurate results when compared to the solutions given directly by the MRCM. Therefore, our numerical experiments consider the approximation of $(\delta\mathbf{u}^n,\delta p^n)$ given by the solution of problem $\mathcal{P}_{\hat{\mathbf{u}}_{0}}$ (\ref{eq:mpm_tf_3}), that is associated with $\epsilon^0$.

We solve Eq. (\ref{eq:mpm_tf_3}) to approximate  $(\hat{\mathbf{u}}_{0},\delta p_0^n)$ by using the MRCM. Note that the conductivity in this equation is $\kappa^m$, thus we can take advantage of the BFs computed at time $t_m$ to solve only one local boundary value problem with trivial Robin boundary conditions for each subdomain. 
Specifically, considering the MRCM additive decomposition of the local approximations, we only need to solve the non-homogeneous part of the solution given by Eq. (\ref{eq:mrc_non_homogeneous}). 
Then, the global unknowns are given by a linear combination of the precomputed BFs, whose coefficients are obtained by solving an inexpensive interface problem. 

With the computed approximation for $(\hat{\mathbf{u}}_{0},\delta p_0^n)$, and hence, for $(\delta\mathbf{u}_0^n,\delta p_0^n)$, we can determine the pair $(\delta\mathbf{u}^n,\delta p^n)$. Then, we can find the solution of problem $\mathcal{P}_{t_n}$ using Eqs. (\ref{eq:soma_p}) and (\ref{eq:soma_u}).
 Let $\bar{\mathbf{u}}=\mathbf{u}^m+\delta\mathbf{u}^n$ denote the approximation of the velocity field at this stage. The resulting approximation is obtained by downscaling, where $\bar{\mathbf{u}}$ defines fluxes on the interfaces of the domain decomposition, which are used as boundary conditions for the following local problems
\begin{equation}\label{mpm_upd}
\left\{
\begin{array}{rll}
\tilde{\mathbf{u}}_h^i &= - \kappa^n\nabla \tilde{p}_h^i &\mbox{in}\ \Omega_i \\
\nabla \cdot \tilde{\mathbf{u}}_h^i &= q^i &\mbox{in}\ \Omega_i \\
\tilde{\mathbf{u}}_h^i \cdot \check{\mathbf{n}}^i&= \bar{\mathbf{u}}_h^i\cdot\check{\mathbf{n}}^i &\mbox{on}\  \partial\Omega_i\cap\partial\Omega\\
\tilde{\mathbf{u}}_h^i \cdot \check{\mathbf{n}}^i&= \bar{\mathbf{u}}_h^i\cdot\check{\mathbf{n}}^i &\mbox{on}\ \Gamma_{ij}\ \forall j
\end{array} \right. 
\end{equation}
for all $\Omega_i,\ i=1,2,\cdots, N$, where $\tilde{\mathbf{u}}$ is the final approximation for $\mathbf{u}^n$. We remark that these local problems are undefined up to a pressure constant. This indeterminacy is removed by imposing a value for the pressure variable at some point of the computational domain.

\subsection{A modified operator splitting scheme}

The operator splitting scheme for two-phase flows calls for an updated velocity field at times $t_n = n\Delta t_p$, for $n=0,1,\dots$. We propose a modification in the algorithm to incorporate the option of choosing between the MRCM and MPM at the elliptic solution steps. 

In order to explain our modified operator splitting algorithm, let  $\{\Phi^m,\Psi^m\}$ denote the set of basis functions $\{\phi_1, \phi_2, \cdots , \phi_{N_U}\}$ and $\{\psi_1, \psi_2, \cdots , \psi_{N_P}\}$ built by the MRCM to solve the problem $\mathcal{P}_{t_m}$ (\ref{eq:mpm_t0}), associated with time $t_m$. 
We compute $p^0$ and $\mathbf{u}^0$ with the MRCM and store the set of BFs  $\{\Phi^0,\Psi^0\}$, such that we can use these basis functions to solve problems $\mathcal{P}_{t_n}$ (\ref{eq:mpm_tf}), at times $t_1,t_2,\dots$. 

The closer the field $\kappa^m$ is to the field $\kappa^n$, the more accurate is
the approximation provided by the MPM \cite{ALI2020125023}. Since the field $\kappa^n$ takes into account an updated saturation, it can be far from $\kappa^m$ depending on the changes due to the displacement of oil by water interface throughout the domain.  
The difference between $\kappa^n$ and $\kappa^m$ is given by $\epsilon=||\kappa^n-\kappa^m||$. We intend to impose a tolerance for $\epsilon$ values in the modified operator splitting scheme to control the difference between $\kappa^n$ and $\kappa^m$. For this purpose, we may need to update the BFs more than once throughout the simulation. Thus, we propose to separate the elliptic solutions into two cases: the case when the BFs are reused by the MPM and the case when a full update of the BFs is required. The latter is computed directly by the MRCM. 

We start a two-phase flow simulation with the solutions $p^0$ and $\mathbf{u}^0$ computed by the MRCM, and the corresponding set of BFs  $\{\Phi^0,\Psi^0\}$ stored. We use these basis functions to solve problems $\mathcal{P}_{t_n}$ (\ref{eq:mpm_tf}), at times $t_n=t_1,t_2,\cdots, t_{m_1-1}$, where $t_{m_1}$ is the first time such that $\epsilon>\eta$ ($\eta$ is the chosen tolerance). At time $t_{m_1}$ we compute $p^{m_1}$ and $\mathbf{u}^{m_1}$ by the MRCM and store the updated set of BFs  $\{\Psi^{m_1},\Phi^{m_1}\}$. Then, we use these BFs to solve problems $\mathcal{P}_{t_n}$, at times $t_n=t_{m_1+1},t_{m_1+2},\cdots, t_{m_2-1}$, where $t_{m_2}$ is the next time when $\epsilon>\eta$, hence we compute $p^{m_2}$, $\mathbf{u}^{m_2}$ and the updated set of BFs by the MRCM. We repeat this procedure until the final simulation time. The MPM-2P algorithm is summarized in Algorithm \ref{algoritmo::mpm}, where $T_{e}$ denotes the total of elliptic solutions computed.

\begin{algorithm} 
\caption{Solving equations (\ref{Darcy})-(\ref{BL2D}) by the modified operator splitting} \label{algoritmo::mpm}
\begin{algorithmic}[1]
\STATE Given $s^{0}(\mathbf x)$, compute $\kappa^{m_0=0}=\lambda(s^{0}(\mathbf x))$
\STATE Compute $p^0(\mathbf x)$ and $\mathbf u^0(\mathbf x)$ from Eq. (\ref{Darcy}) by using the MRCM
\STATE Store the set of BFs $\{\Psi^{m_0=0},\Phi^{m_0=0}\}$
\STATE Set $n=1$, $\ell=0$, and $\epsilon=\eta$
\WHILE {$n < T_e$ }
	\FOR{$k\in\{1,\cdots,C_{n-1}\}$ }
		\STATE $t_{n-1,k} =  t_{n-1} + k\Delta t_s$ 
		\STATE Solve Eq. (\ref{BL2D}) to compute $s(\mathbf x, t_{n-1,k})$ 
	\ENDFOR
	\STATE Given $s^{n}(\mathbf x)$, update $\kappa^n=\lambda(s^{n}(\mathbf x))$
    \IF{$\epsilon>\eta$} 
    		\STATE $\ell=\ell+1$
    		\STATE Compute $p^{n=m_{\ell}}(\mathbf x)$ and $\mathbf u^{n=m_{\ell}}(\mathbf x)$ from Eq. (\ref{Darcy}) by using the MRCM
    		\STATE Store the updated set of BFs $\{\Psi^{m_{\ell}},\Phi^{m_{\ell}}\}$
    		\STATE Update $\kappa^{m_{\ell}}=\lambda(s^{n=m_{\ell}}(\mathbf x))$
    	\ELSE
    		\STATE Compute $p^n(\mathbf x)$ and $\mathbf u^n(\mathbf x)$ from Eq. (\ref{Darcy}) with the MPM, reusing BFs $\{\Psi^{m_{\ell}},\Phi^{m_{\ell}}\}$
    	\ENDIF
    \STATE Compute $\epsilon=\parallel\kappa^{n}-\kappa^{m_{\ell}}\parallel$
  \STATE $n =  n +1$
 \ENDWHILE
\end{algorithmic}
\end{algorithm}

\subsection{Computational cost of the MPM-2P}\label{subsection:cost}

To compare the computational cost of the MRCM and MPM-2P in the solution of the elliptic equations arising within the operator splitting algorithm we start by computing the number of BFs required by them, considering a problem with a domain decomposition with $N=N_x\times N_y$ subdomains (2D).
In order to find the number of BFs required for each method in the approximation of the two-phase flow problem, let us consider that a total of $T_{e}$ elliptic solutions need to be computed. Note that:

\begin{itemize}

\item \textit{Number of BFs required for an elliptic solution: }

The number of BFs required by the MRCM in each subdomain is $4\times (k_U+k_P)$ homogeneous BFs for each one of the edges of the subdomain, plus one non-homogeneous basis function. This number may be different for distinct subdomains due to local choices of degrees of freedom per interface. Let $\hat{N}$ be the total number of homogeneous BFs required by the MRCM, and hence, the total amount of BFs computed by the MRCM is $\hat{N} + N$ (in a serial mode implementation).

The MPM requires only the calculation of the basis function for the non-homogeneous part of the solution in each subdomain. Therefore, we have a total of $N$ BFs.

\item  \textit{Number of BFs for the coupled flow and transport problem:}

The number of BFs required by the MRCM for two-phase flows is $(\hat{N} + N)\times T_{e}$. To compute the total of BFs required by the MPM-2P we have to separate the cases when the basis functions are reused from the cases when a full update is required. 

\begin{enumerate}

\item Let  $T_m$ be the total number of updates required by the MPM-2P (associated with the counter $\ell$ at line 12 of Algorithm \ref{algoritmo::mpm}). If we compute each update with the MRCM (considering the same number of BFs), the total number of BFs required by the updates of the MPM-2P is $(\hat{N}+N)\times T_m$. 

\item The total of BFs computed when reusing the basis functions is $N\times(T_{e}-T_m)$.

\end{enumerate}
Therefore, the total number of BFs computed by the MPM-2P is $(\hat{N}+N)\times T_m+N\times(T_e-T_m)$.

\end{itemize}

To estimate the overall cost of the methods we have to consider the cost of computing the BFs, downscaling, and a global interface problem. Let $\mathcal{C}_{BF}$, $\mathcal{C}_{DS}$ and $\mathcal{C}_{I}$ be, respectively, the estimated computational cost to compute one basis function, the downscaling in a subdomain, and the global interface problem. We define the cost estimate of the MRCM as follows:
 \begin{equation}
\begin{array}{rll}
\text{cost(MRCM)}&=\left[\mathcal{C}_{BF}\times(\hat{N}+N) + \mathcal{C}_{DS} \times N + \mathcal{C}_{I}\right]\times T_{e}\vspace{0.2cm}\\ 
&\approx\mathcal{C}_{BF}\times(\hat{N}+2\times N) \times T_e.
\end{array}
\end{equation}
This approximation follows from the fact that the computational cost of the interface problem is typically negligible when compared to the cost of computing BFs \cite{Paolathesis,PaolaJCP}. Furthermore, the downscaling step has essentially the same cost of computing one basis function at each subdomain ($\mathcal{C}_{DS}\approx\mathcal{C}_{BF}$).
Thus, the cost estimate of the MPM-2P is given by:
 \begin{equation}
\begin{array}{rll}
\text{cost(MPM-2P)}&=\left[\mathcal{C}_{BF}\times(\hat{N}+N) + \mathcal{C}_{DS}\times N\right]\times T_m\vspace{0.2cm}\\
&\qquad+(\mathcal{C}_{BF} +\mathcal{C}_{DS})\times N\times(T_{e}-T_m)+\mathcal{C}_{I}\times T_{e}\vspace{0.2cm}\\ 
&\approx \left[\mathcal{C}_{BF}\times(\hat{N}+2\times N)\right] \times T_m+2\times\mathcal{C}_{BF}\times N\times(T_{e}-T_m).
\end{array}
\end{equation}

We define a quantity to indicate the relation between the computational cost of the methods. The following quantity measures the Relative Cost Reduction (RCR) accomplished by the MPM-2P when compared with the approximation of two-phase flows directly by the MRCM.
 \begin{equation}
\begin{array}{rl}
\text{RCR}&=\dfrac{\text{cost(MRCM)-cost(MPM-2P)}}{\text{cost(MRCM)}} 100\% \vspace{0.4cm}\\ 
&=\dfrac{\left[\mathcal{C}_{BF}\times(\hat{N}+2\times N)\right] \times (T_e - T_m) - 2\times\mathcal{C}_{BF}\times N\times(T_{e}-T_m)}{\mathcal{C}_{BF}\times(\hat{N}+2\times N) \times T_e} 100\%
\vspace{0.4cm}\\ 
&=\dfrac{T_{e}-T_m}{T_{e}}\left[1-\dfrac{2\times N}{\hat{N}+2\times N}\right]100\%.
\end{array}
\end{equation}

Let us consider as an example a domain decomposition of $4\times 4$ subdomains that is used in some of the numerical experiments below.
If we consider the MRCM with constant interface spaces for both flux and pressure, i.e. $k_U = k_P = 1$, the total number of homogeneous BFs to be computed is $\hat{N}=96$ (considering the physical boundary conditions). In order to find the RCR for a two-phase flow problem, let us consider that a total of $T_e=3500$ elliptic solutions need to be computed. This is typically the order of the number of elliptic solutions needed to reach water breakthrough in some of our simulations. We find in our numerical experiments that usually, less than 10 updates are required by the MPM-2P for this type of problem. Therefore, the RCR is given by
 \begin{equation}
\begin{array}{rl}
\text{RCR}&=\dfrac{T_{e}-T_m}{T_{e}}\left[1-\dfrac{2\times N}{\hat{N}+2\times N}\right]100\%\vspace*{0.4cm}\\
&=\dfrac{3500-10}{3500}\left[1-\dfrac{2\times 16}{96+2\times 16}\right]100\%\approx 74.79\%
\end{array}
\end{equation}
The values attained by the cost function for each one of our numerical experiments are shown in the following section. We find that the MPM-2P presents outstanding speed-up. It reduces significantly the cost of the simulation of two-phase flows when compared to the traditional operator splitting combined with the MRCM. 
A RCR of $68.60\%$ is the least value that we find in our numerical experiments. The more basis functions we consider the greater is the advantage of using the MPM-2P. 
 
\section{Numerical Results} \label{Sec:4}

In this section, we present numerical simulations to investigate the accuracy as well as the computational cost of the MPM-2P.
We consider challenging two-phase flow problems, with high-contrast permeability fields and water-oil finger growth in a homogeneous medium. 

In all simulations, we set the relative permeability curves $k_{ro}=(1-s)^2$ and $k_{rw}=s^2$, and hence, the fractional flow function is given by
\begin{equation}
f(s) = \dfrac{M s^2}{M s^2 + (1-s)^2},
\end{equation}
where $ M = {\mu_o}/{\mu_w}$. The time is expressed in PVI (Pore Volume Injected) \cite{zchen_book}, and the results are presented in terms of the number of elliptic solutions. The downscaling procedure used to compute a conservative solution for the MRCM approximation is the Stitch method presented in \cite{guiraldello2019downscaling}.

\subsection{A Gaussian permeability field}\label{subsec:A Gaussian permeability field}

\begin{figure}[htbp]
\centering
\includegraphics[scale=0.6]{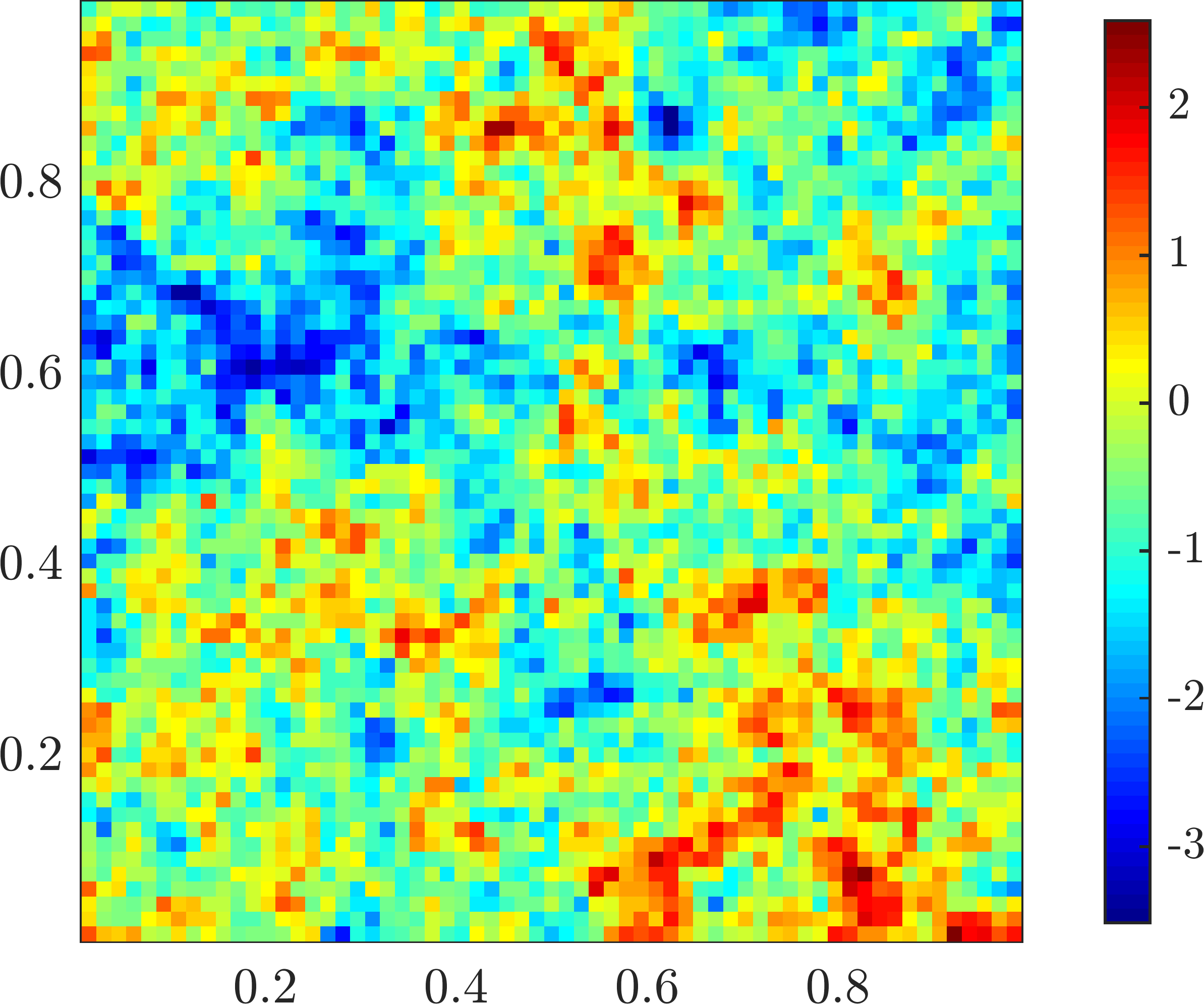}
\caption{Gaussian permeability field (log-scaled) with a permeability contrast of $K_{\max}/K_{\min}\approx10^6$.}
\label{K_slab}
\end{figure}

In the first example, we will consider a slab geometry problem with a Gaussian permeability field. Our initial assumption is that the reservoir is fully saturated with oil. Water is then injected at a constant rate. Moreover, here $M = 40$ in the definition of $f(s)$. The computational domain is taken as a square $[0,1]\times[0,1]$ containing $64\times 64$ fine grid cells. There is a Dirichlet boundary on the left ($p=1$) and right ($p=0$). The top and bottom are no flow (Neumann) boundary conditions. The domain is divided into $4\times 4$ subdomains with each subdomain having $16\times 16$ fine cells. There are no source terms taken into account for this example.
The permeability data is considered to be $K(\mathbf{x}) = 0.8\ e^{\delta \xi(\mathbf{x})}$, where $\delta=2.5$ for a permeability contrast of $K_{\max}/K_{\min}\approx10^3$ and $\delta=4.5$ for $K_{\max}/K_{\min}\approx10^6$. The field $\xi(\mathbf{x})$ is a self similar Gaussian distribution having zero mean and the covariance function
given by $C(\mathbf{x},\mathbf{y})=|\mathbf{x}-\mathbf{y}|^{-1/2}$. A sample permeability field is shown in Fig. \ref{K_slab}.

We will discuss the relative error obtained by the MPM-2P and the MRCM for both the flux and the saturation with respect to a reference fine grid solution. The updates of the BFs for the MPM-2P consider the same set-up of the approximation given directly by the MRCM, with a tolerance of $\eta = 10^{-2}$ for the values of $\epsilon$, in line with \cite{ALI2020125023}. The errors are shown as a function of the number of elliptic solutions. The flux error is computed in terms of the $L^2(\Omega)$ norm while the saturation is in terms of the $L^1(\Omega)$ norm. 
In addition to that, we will discuss the saturation profiles obtained at the breakthrough time.

In this example, the interface spaces for the MRCM are the simplest possible, being constant for both pressure and flux. Additionally, we use an intermediate length scale $h\leq \bar{H}\leq H$ to define the constant polynomials at the interfaces of the subdomains. We test in our numerical experiments two choices: $\bar{H}=H=16h$, that is the classic choice of one constant basis function per subdomain interface, and  $\bar{H}=H/2=8h$, that represents a division of each subdomain interface into two parts, each one containing a constant basis function. The MRCM solution with constant interface spaces along with the algorithmic function set as $\alpha(\mathbf{x}) = 1$ is equivalent to the solution yielded by the MuMM \cite{pereira}.

\begin{figure}[htbp]
\centering
\includegraphics[scale=0.46]{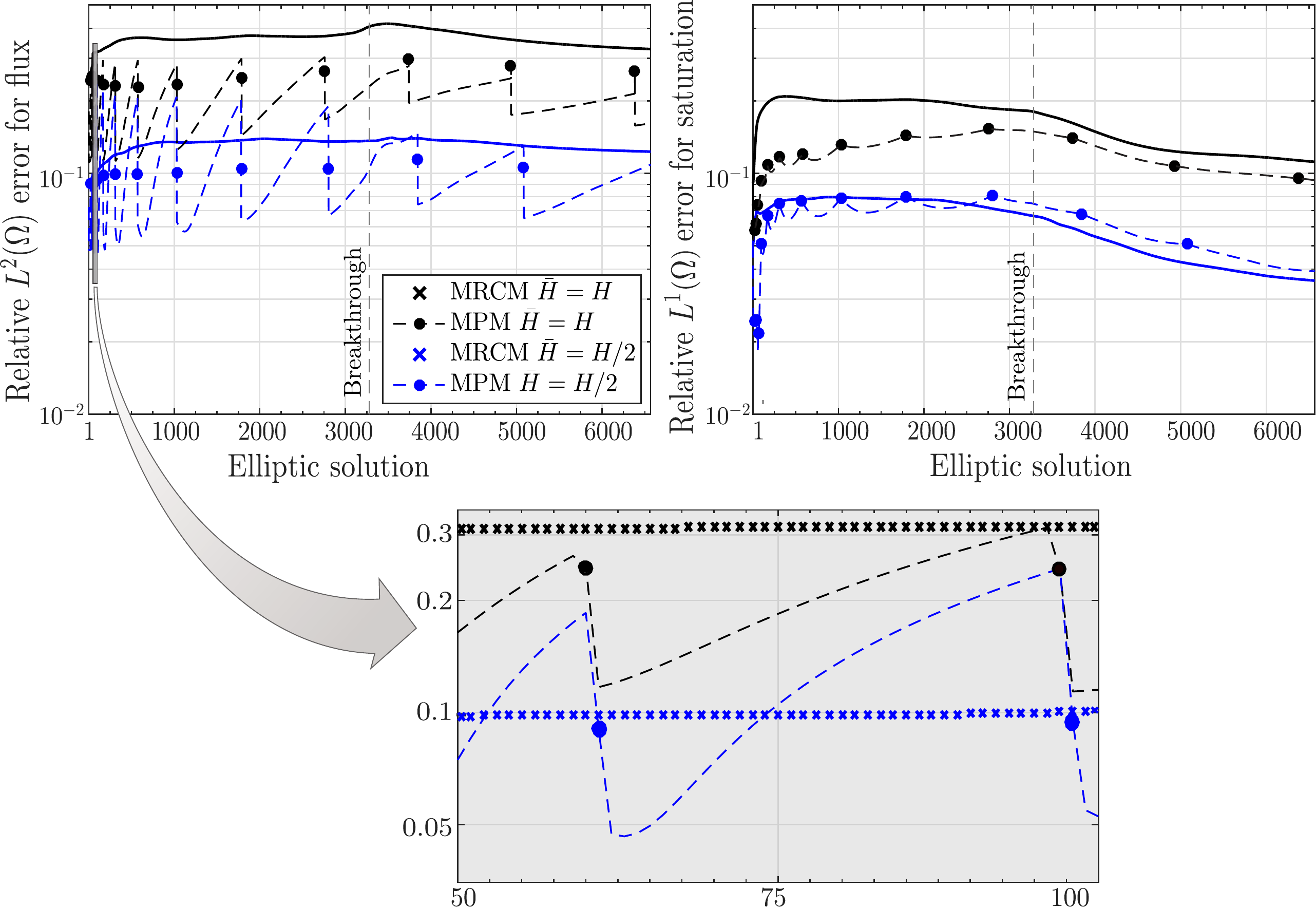}
\caption{Relative errors for the Gaussian permeability field with contrast of $K_{\max}/K_{\min} \approx 10^3$. Top: the relative $L^2(\Omega)$ error norm for the flux (left) and the relative $L^1(\Omega)$ error norm for the saturation (right). Bottom: the zoomed version of the flux error shown in the top left. We consider $\bar{H}=H$ and $\bar{H}=H/2$ and compare the MRCM and the MPM-2P. The nodes on each of the MPM-2P curves indicate the times when the BFs are updated. The breakthrough time is illustrated by a vertical dashed line. Note that the MPM-2P is significantly less expensive than the MRCM.}
\label{K_3_Cn_1}
\end{figure}

Figure \ref{K_3_Cn_1} indicates the relative error for the flux and saturation obtained for the permeability contrast of $10^3$ and setting $C_n = 1,\ \forall n$, which means that the elliptic solution has to be computed again after each transport step. There is a comparison between the results obtained by using the MRCM and the MPM-2P, both considering $\bar{H}=H$ and $\bar{H}=H/2$. 
The MPM-2P with $\bar{H}=H$ (and $\bar{H}=H/2$ respectively) is more accurate than the MRCM with $\bar{H}=H$ (and $\bar{H}=H/2$ respectively). The lowest accuracy is observed for the MRCM approximation with $\bar{H} = H$. Moreover, the highest accuracy can be observed in the case of the MPM-2P with $\bar{H} = H/2$. 
Remember that, in the case of the MRCM, the set of BFs get recomputed at every elliptic update, while, in the case of the MPM-2P algorithm, the BFs are recomputed only when  $\epsilon > 10^{-2}$. The nodes on the curves corresponding to the relative error obtained using the MPM-2P indicate these elliptic updates. The zoomed version in Fig. \ref{K_3_Cn_1} indicates precisely the advantages of using the MPM-2P. Each blue or black cross in the MRCM indicates the recalculation of the BFs while in the case of MPM-2P, the two black nodes, and the two blue nodes indicate elliptic updates where the BFs were recomputed. With respect to the breakthrough time, of the fine grid solution, the MRCM computes the set of BFs $3283$ times while the MPM-2P computes the set of BFs $10$ times (the initial set plus 9 updates). This is where the real computational advantage of the MPM-2P is observed.

\begin{figure}[htbp]
\centering
\includegraphics[scale=0.75]{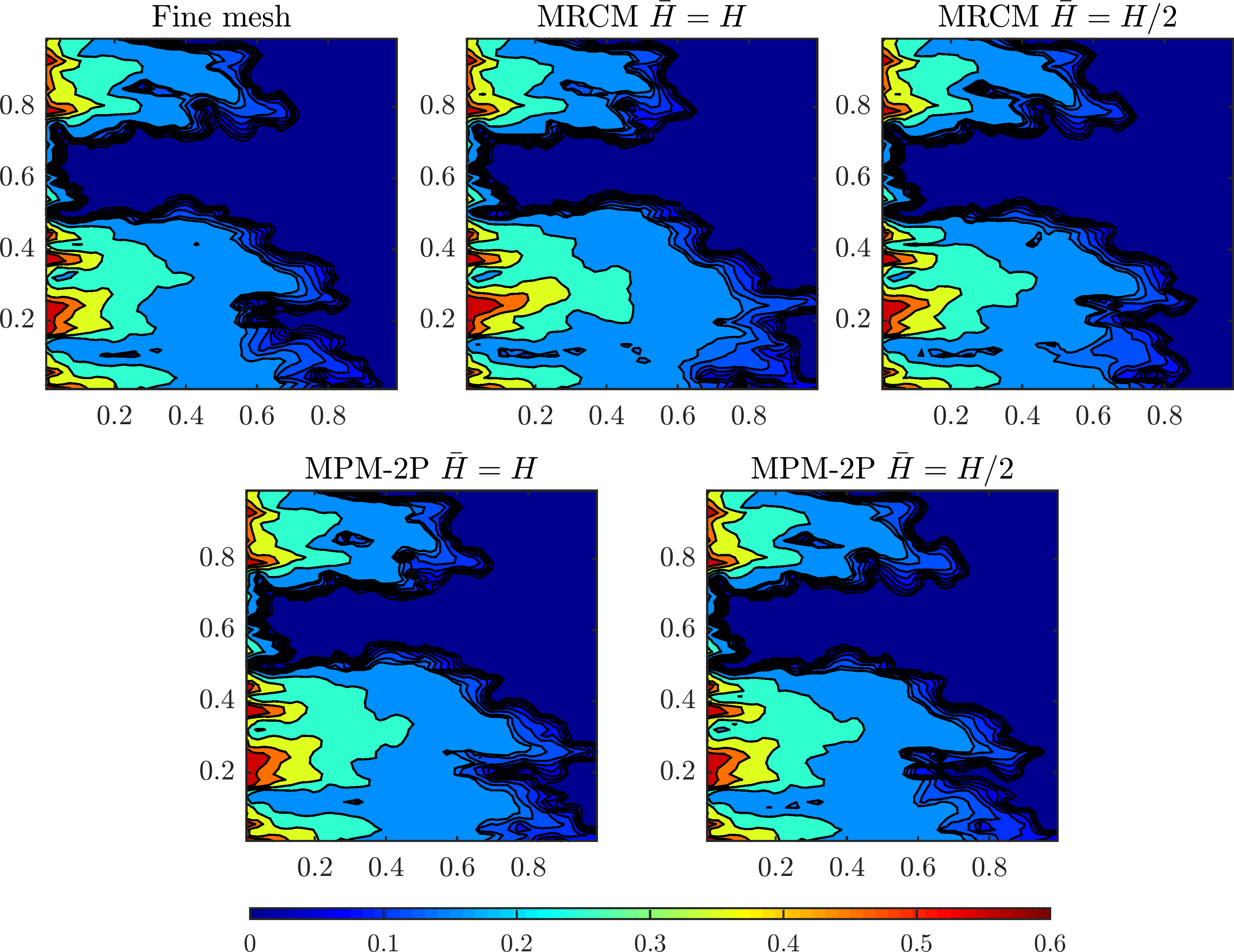}
\caption{Saturation profiles at breakthrough time $T_{\text{PVI}}=0.12$ (that corresponds to the elliptic time step number $3283$) for the Gaussian permeability field with contrast of $K_{\max}/K_{\min}\approx10^3$. First line, left to right: fine grid solution; MRCM with $\bar{H}=H$; MRCM with $\bar{H}=H/2$. Second line shows the solutions for the MPM-2P, with $\bar{H}=H$ (left) and $\bar{H}=H/2$ (right).}
\label{K_3_Sat}
\end{figure}

Figure \ref{K_3_Sat} shows the saturation profiles for the previous experiment obtained at the breakthrough time $T_{\text{PVI}}=0.12$ (elliptic time step number $3283$). We can compare how close or accurate the approximation obtained by using MPM-2P and the MRCM with different $\bar{H}$ values are to the fine scale solution. We can observe that the approximations that consider $\bar{H} = H/2$ capture the saturation profile in more accurate manner compared to the approximations with $\bar{H}=H$. Therefore, the MPM-2P approximation with $\bar{H}=H/2$ helps to lower the computational cost as well as produce approximations which are accurate.

\begin{figure}[htbp]
\centering
\includegraphics[scale=0.45]{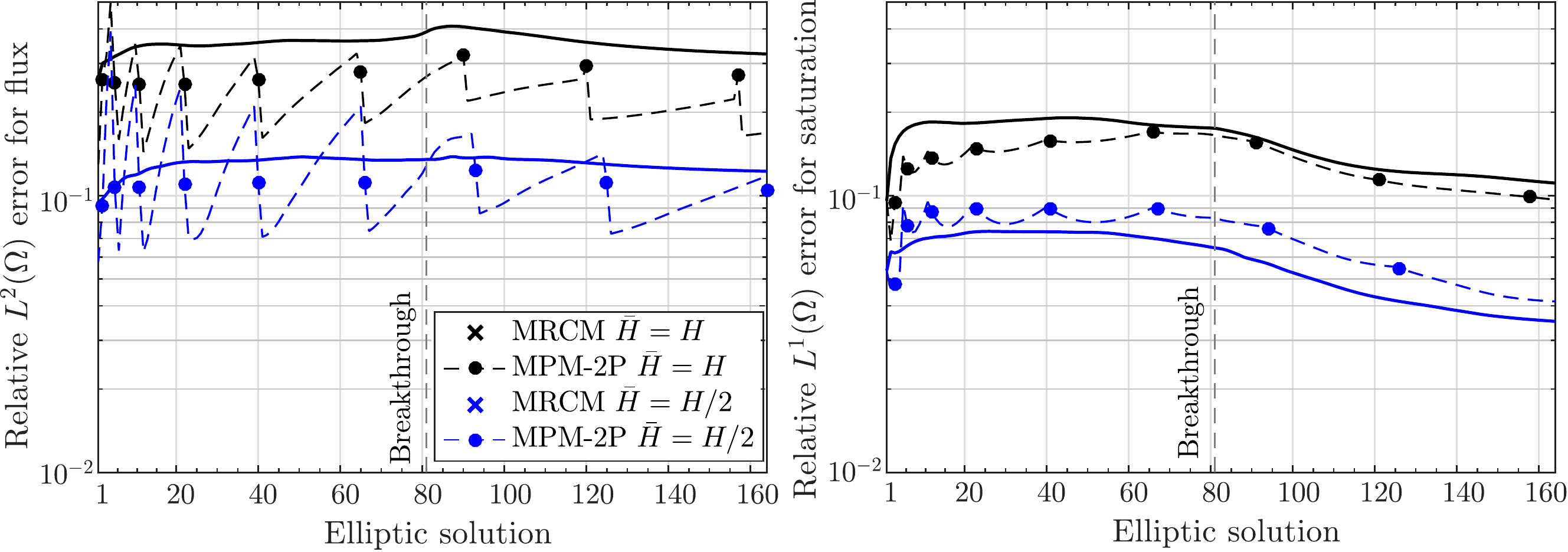}
\caption{Relative errors for the Gaussian permeability field with contrast of $K_{\max}/K_{\min} \approx 10^3$ and $C_n \approx 20,\  \forall n$. Relative $L^2(\Omega)$ error for  flux (left) and $L^1(\Omega)$ error for saturation (right). We consider $\bar{H}=H$ and $\bar{H}=H/2$ and compare the MRCM and the MPM-2P. The nodes on each of the MPM-2P curves indicate the times when the BFs are updated.  The results are essentially the same as those attained with $C_n=1,\  \forall n$.}
\label{K_3_Cn_20}
\end{figure}

We perform the same study by considering approximately 20 transport steps between successive elliptic updates ($C_n \approx 20,\  \forall n$), which means that a smaller number of elliptic solutions will be required during the simulation. 
Figure \ref{K_3_Cn_20} shows a comparison between the results obtained by using the MRCM and the MPM-2P. Results are similar to the case with $C_n=1,\  \forall n$, where we note a higher accuracy for the MPM-2P with $\bar{H} = H/2$ as compared to the other cases. 

\begin{figure}[htbp]
\centering
\includegraphics[scale=0.45]{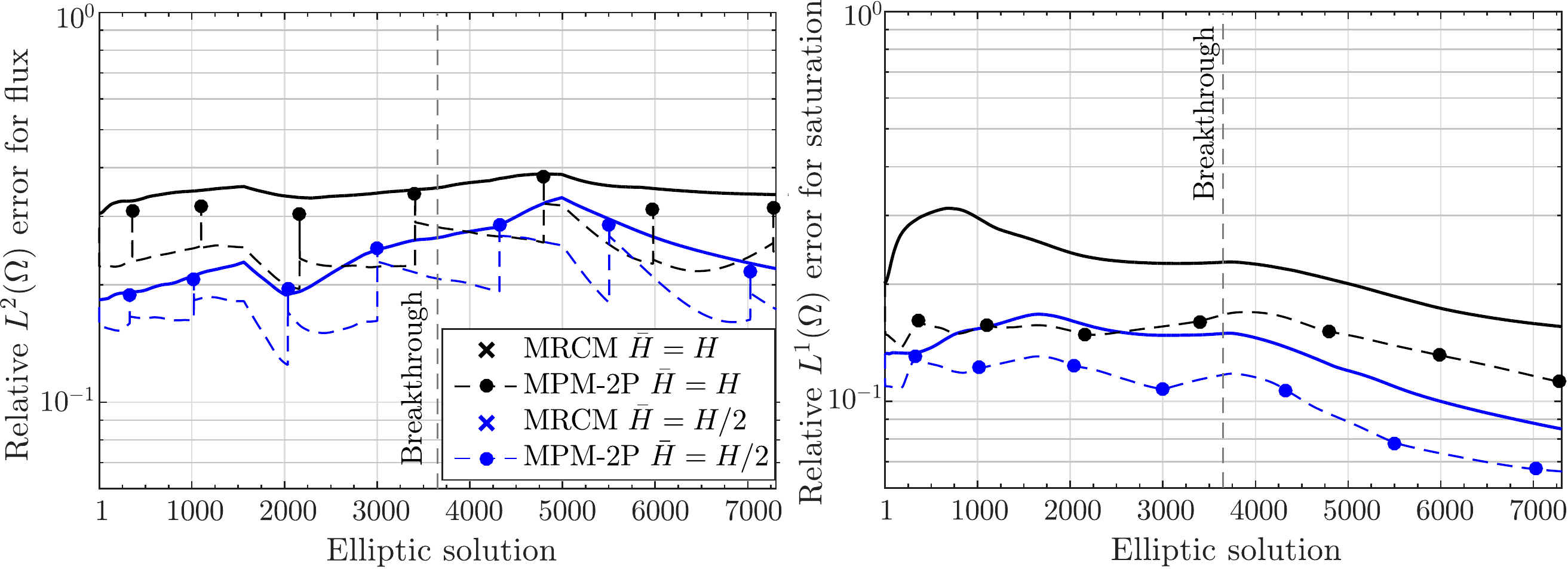}
\caption{Relative errors for the Gaussian permeability field with contrast of $K_{\max}/K_{\min} \approx 10^6$ and $C_n =1,\ \forall n$. Relative $L^2(\Omega)$ error for flux (left) and $L^1(\Omega)$ error for saturation (right). We consider $\bar{H}=H$ and $\bar{H}=H/2$ and compare the MRCM and the MPM-2P. The nodes on each of the MPM-2P curves indicate the times when the BFs are updated.}
\label{K_6_Cn_1}
\end{figure}

\begin{figure}[htbp]
\centering
\includegraphics[scale=0.75]{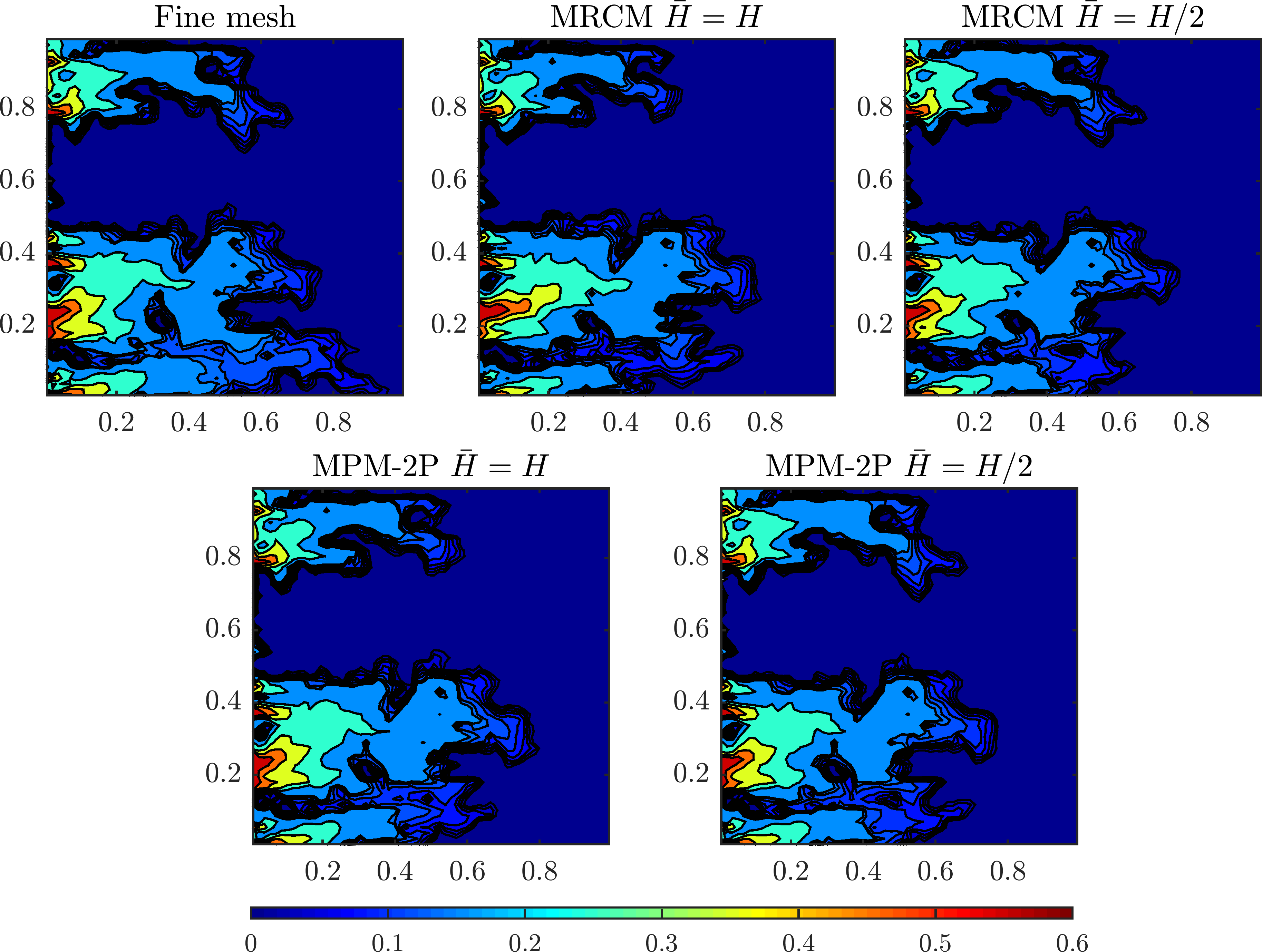}
\caption{Saturation profiles at breakthrough time $T_{\text{PVI}}=0.10$ (that corresponds to the elliptic time step number $3652$) for the Gaussian permeability field with contrast of $K_{\max}/K_{\min}\approx10^6$. First line, left to right: fine grid solution; MRCM with $\bar{H}=H$; MRCM with $\bar{H}=H/2$. Second line shows the solutions for the MPM-2P, with $\bar{H}=H$ (left) and $\bar{H}=H/2$ (right).}
\label{K_6_Sat}
\end{figure}

Next, we conducted similar experiments with a higher contrast permeability ratio where $K_{\max}/K_{\min} \approx 10^{6}$. With the new permeability field, the problem becomes more difficult and presents a challenge for the multiscale methods. Results obtained for this case with $C_n = 1,\ \forall n$, can be seen in Fig. \ref{K_6_Cn_1}.  We observe similar patterns in the relative error for the flux and saturation as discussed for the case where $K_{\max}/K_{\min} \approx 10^{3}$. 
The MPM-2P with $\bar{H} = H/2$ gives the best approximation in terms of accuracy and computational cost. The saturation profiles at the breakthrough time $T_{\text{PVI}}=0.10$ (elliptic time step number $3652$) are shown in Fig. \ref{K_6_Sat}. Here also we can make a similar observation as we did for Fig. \ref{K_3_Sat}.

\begin{figure}[htbp]
\centering
\includegraphics[scale=0.45]{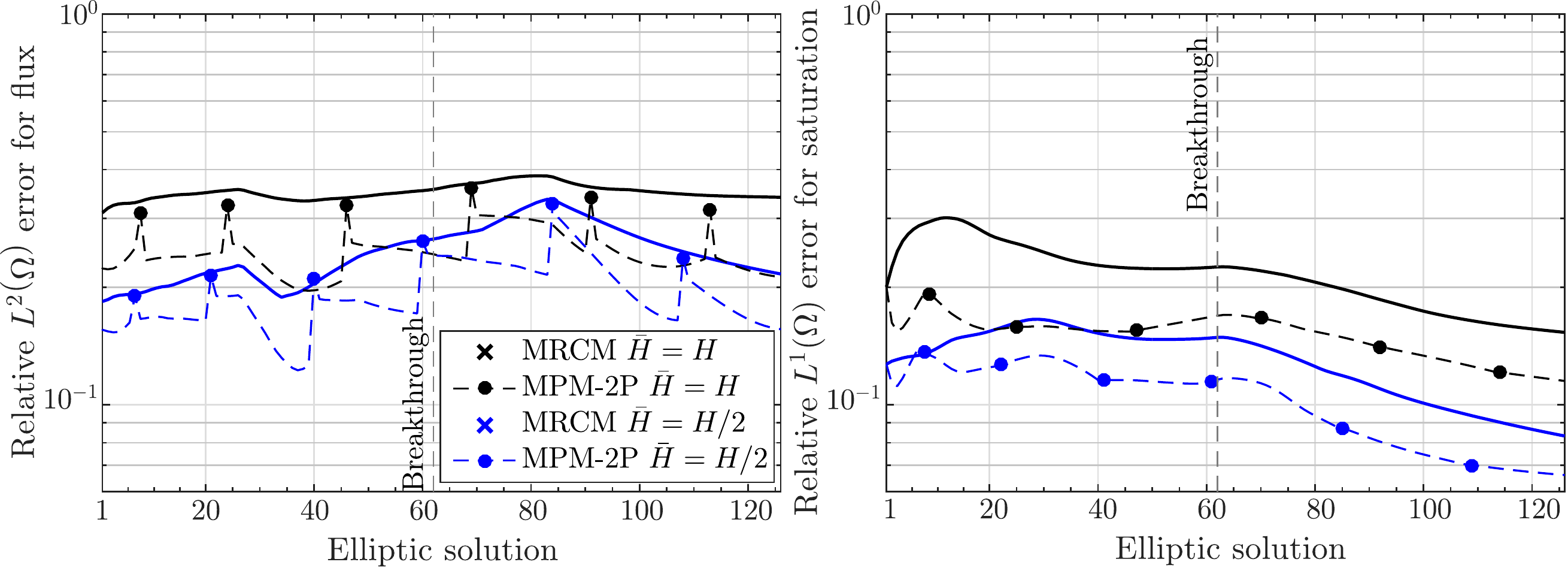}
\caption{Relative errors for the Gaussian permeability field with contrast of $K_{\max}/K_{\min} \approx 10^6$ and $C_n \approx 20,\  \forall n$. Relative $L^2(\Omega)$ error for flux (left) and $L^1(\Omega)$ error for saturation (right). We consider $\bar{H}=H$ and $\bar{H}=H/2$ and compare the MRCM and the MPM-2P. The nodes on each of the MPM-2P curves indicate the times when the BFs are updated. The MPM-2P with $\bar{H} = H/2$ is a good balance between accuracy and computational cost.}
\label{K_6_Cn_20}
\end{figure}

We also perform the study for the higher permeability contrast by considering approximately 20 transport steps between successive elliptic updates ($C_n \approx 20,\ \forall n$), which are summarized in Fig. \ref{K_6_Cn_20}. The obtained results are similar to the case with $C_n=1,\ \forall n$, where we can conclude that the MPM-2P with $\bar{H} = H/2$ is a good balance between accuracy and computational cost.

Table \ref{RCR_Gauss} presents the Relative Cost Reduction obtained by the MPM-2P for the previous experiments. We can observe that with the use of the MPM-2P we are able to obtain an exceptional speed-up: ranging from around 68.60\% to 85.60\% for both the types of permeability ratios. Hence, we can make one more conclusion that the MPM-2P can significantly reduce the computational cost of solving the two-phase flow problem irrespective of the contrast of the Gaussian permeability field considered. This is a noteworthy achievement when considering its application to solving the two-phase problem with a  more challenging permeability data field.


 \begin{table}[h!]
     \centering
	\caption{Relative Cost Reduction obtained by the MPM-2P for the experiments with the Gaussian permeability field.}
		\label{RCR_Gauss}
    	\resizebox{0.56\textwidth}{!}{\begin{minipage}{3.25in}
     		\begin{tabular}{ c c c c c c c c c } 
		\hline
			& \multicolumn{4}{c}{\ $K_{\max}/K_{\min} = 10^3$\qquad\qquad} & \multicolumn{4}{c}{\ $K_{\max}/K_{\min} = 10^6$}\\
      		 & \multicolumn{2}{c}{$C_n = 1$\qquad}  & \multicolumn{2}{c}{$C_n \approx 20$}& \multicolumn{2}{c}{$C_n = 1$\qquad}  & \multicolumn{2}{c}{$C_n \approx 20$}\\
      		 \hline
       		$\bar{H}=H$ & \multicolumn{2}{c}{74.77\%} &\multicolumn{2}{c}{68.60\%} & \multicolumn{2}{c}{74.90\%} &\multicolumn{2}{c}{69.05\%}\\
       		$\bar{H}=H/2$ & \multicolumn{2}{c}{85.45\%} &\multicolumn{2}{c}{78.40\%} & \multicolumn{2}{c}{85.60\%} &\multicolumn{2}{c}{78.91\%}\\     		
      		 \hline
    		 \end{tabular}
      \end{minipage}}
 \end{table}    


\subsection{A high-contrast permeability field}\label{subsec:A high-contrast permeability field}

\begin{figure}[htbp]
    	\centering
	\includegraphics[scale=0.375]{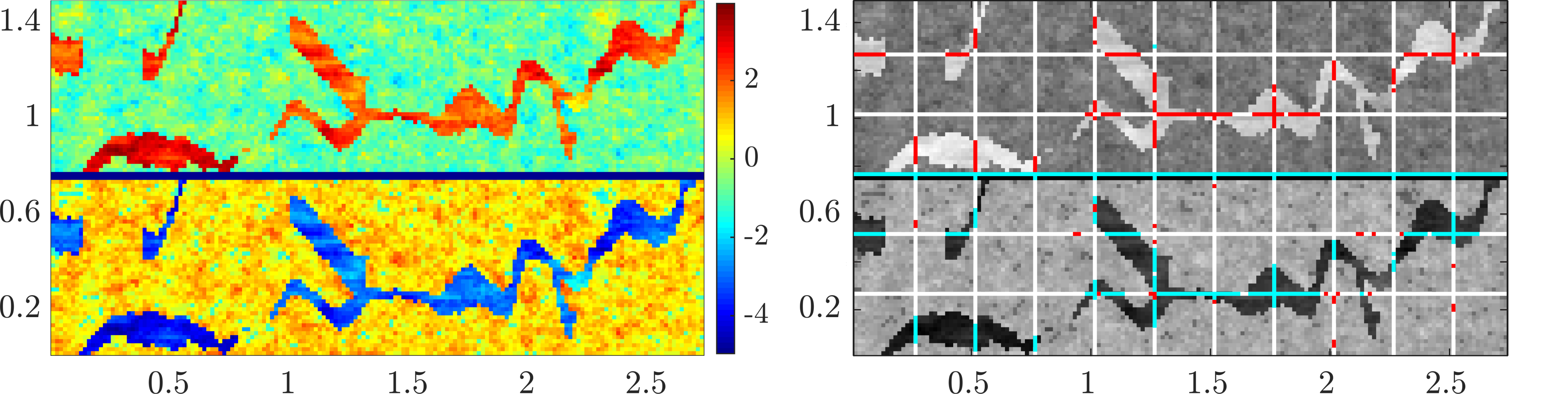}
	    \caption{Log-scaled permeability field (left) and map of the permeability variations at the boundaries of the subdomains (right). On the right, the red color identifies regions of high-permeability and the cyan color represents the regions of low-permeability. Note that the channelized structures and barriers are well captured by this procedure.}
	 \label{fig:MPM_1}
\end{figure}

The second experiment considers a high-contrast permeability field containing a high-permeable channel and a low-permeable region, as illustrated in Fig. \ref{fig:MPM_1} (left). This permeability field is a modification of one of the layers of the SPE-10 project \cite{christie2001tenth}, built to benchmark the methods developed in \cite{bifasico}. This field is very challenging to multiscale methods, since it combines both channels of high permeability and barriers of low permeability in the same problem. The domain $\Omega=[0,33/12]\times[0,3/2]$ is divided into $11\times6$ subdomains with $15\times15$ cells into each one. The flow is established by imposing unit flow at the left boundary and zero pressure at the right boundary along with no-flow at top and bottom. No source terms are considered. Here, we also consider that the porous medium is initially filled with oil and water is injected at a constant rate. The viscosity ratio is set to be $M=40$.

Our objective is to compare the approximations provided by the MPM-2P and those obtained purely by the MRCM. To solve this difficult problem, we consider an improved version of the MRCM, the adaptive MRCM, as presented originally in \cite{bifasico}. The adaptive version of the MRCM automatically sets the parameter $\alpha(\mathbf x)$ on the interfaces of the domain decomposition, depending if the interfaces are crossed by heterogeneities such as high permeable channels or low permeable barriers. The values of $\alpha(\mathbf x)$ are set based on a threshold function, which is illustrated in Fig. \ref{fig:MPM_1} (right) for the chosen permeability field. In the MPM-2P method, the recalculation of the basis functions, when needed, will be performed by the same adaptive MRCM. Both methods (MPM-2P and $a$MRCM) will be tested with two different interface spaces, a linear polynomial interface space (denoted by the suffix -POL), and the interface spaces based on physics, as presented in \cite{rocha2020interface} (denoted by the suffix -PBS). We remark that the cost and accuracy of the linear polynomial interface space are about the same as using constant spaces with $\bar{H} = H/2$, as performed in the previous section. The interested reader is referred to \cite{bifasico,rocha2020interface, rochaenhanced} for more details about the parameters of adaptive MRCM and about the construction of the interface spaces based on physics.

%

\begin{figure}[htbp]
    	\centering
	\includegraphics[scale=0.45]{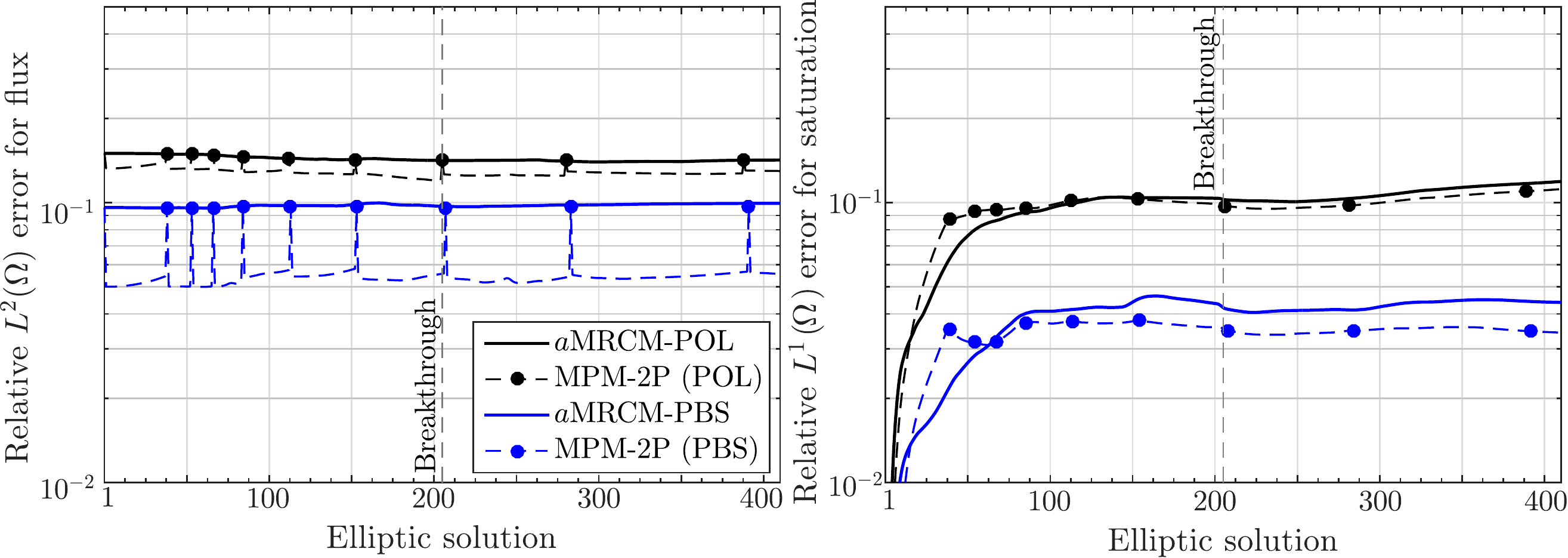}
	    \caption{Relative $L^2(\Omega)$ errors for flux (left) and relative $L^1(\Omega)$ errors for saturation (right). We compare both MPM-2P (with $\eta=10^{-2}$) and $a$MRCM, for both choices of interface spaces (POL and PBS). The nodes on each of the MPM-2P curves indicate the times when the BFs are updated. Note the improved accuracy of both methods when using physics-based interface spaces.}
    \label{fig:MPM_2}
\end{figure}

As in the previous results, we have chosen $C_n \approx 20$, which means nearly $20$ transport steps between successive elliptic updates for all methods considered. 
Figure \ref{fig:MPM_2} shows the relative $L^2(\Omega)$ errors for flux (left) and relative $L^1(\Omega)$ errors for saturation (right) as a function of the number of elliptic solutions performed. The errors are computed with respect to a reference fine grid solution. In this figure, the breakthrough time for the fine grid solution is indicated by a vertical dashed line.

Remember that the BFs are fully rebuild every elliptic solution when using purely the adaptive MRCM. When using the MPM-2P, however, these updates are drastically reduced, being required only when $\epsilon > \eta$. In such cases, the same adaptive MRCM is used to update the BFs.  We compare the errors of the MPM-2P (with $\eta=10^{-2}$) and the adaptive MRCM with the fine (undecomposed) grid solution, that are displayed in Fig. \ref{fig:MPM_2}. The nodes appearing on each of the MPM-2P curves indicate the times when the BFs are updated. One can see in these results that the set of BFs was computed 10 times by the MPM-2P (the initial set plus 9 updates), that is significantly less than the total of $412$ full updates required by the adaptive MRCM.
As expected, the $a$MRCM-PBS (and respectively the MPM-2P that uses the $a$MRCM-PBS) is more accurate than the $a$MRCM-POL (respectively the MPM-2P using the $a$MRCM-POL), but most importantly, the MPM-2P yield  results that are more accurate than the adaptive MRCM alone.
Note that the flux errors of the MPM-2P tend to be slightly lower than those produced by the $a$MRCM, with a subtle increase when the $a$MRCM is invoked to update the basis functions. The error rapidly drops after every BFs full recalculations. By these results, one can  see that the MPM-2P also benefits from the physics-based spaces.

\begin{figure}[htbp]
    	\centering
	\includegraphics[scale=0.45]{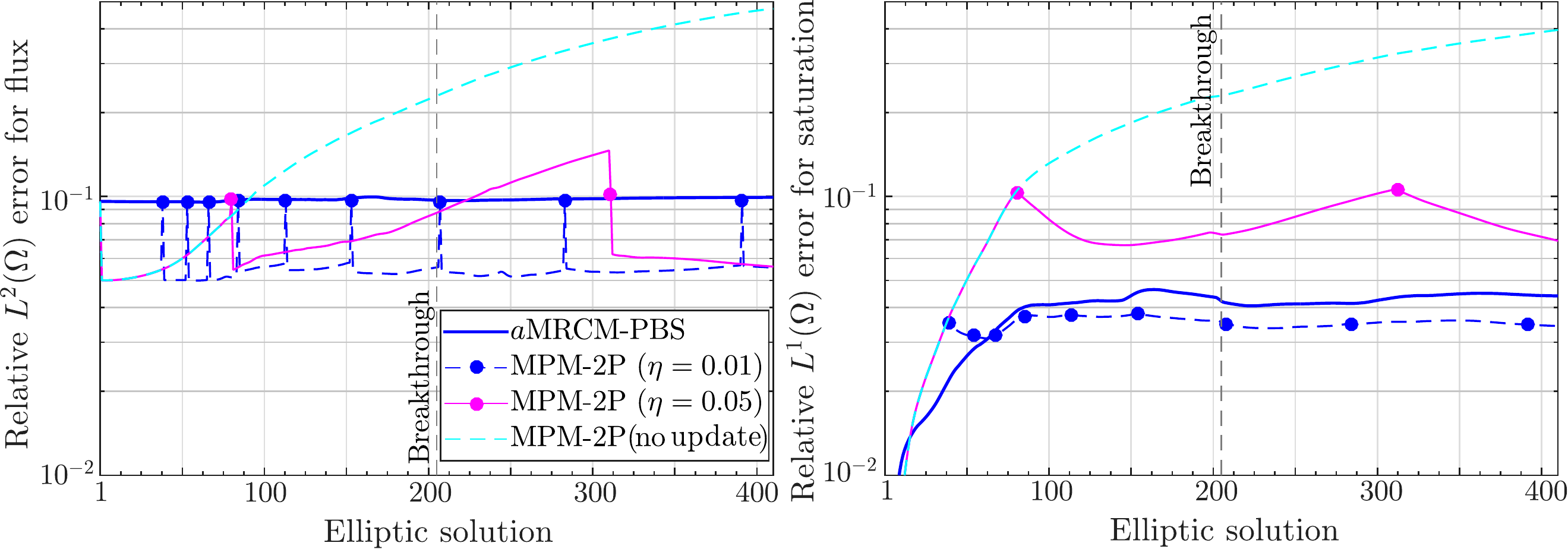}
	    \caption{Relative $L^2(\Omega)$ errors for flux (left) and relative $L^1(\Omega)$ errors for saturation (right). We compare the $a$MRCM-PBS with the MPM-2P with different tolerances: with no BFs updates, $\eta=0.01$, and $\eta=0.05$. The nodes on each of the MPM-2P curves indicate the times when the BFs are updated.}
	\label{fig:MPM_4}
\end{figure}

\begin{figure}[htbp]
    	\centering
	\includegraphics[scale=0.68]{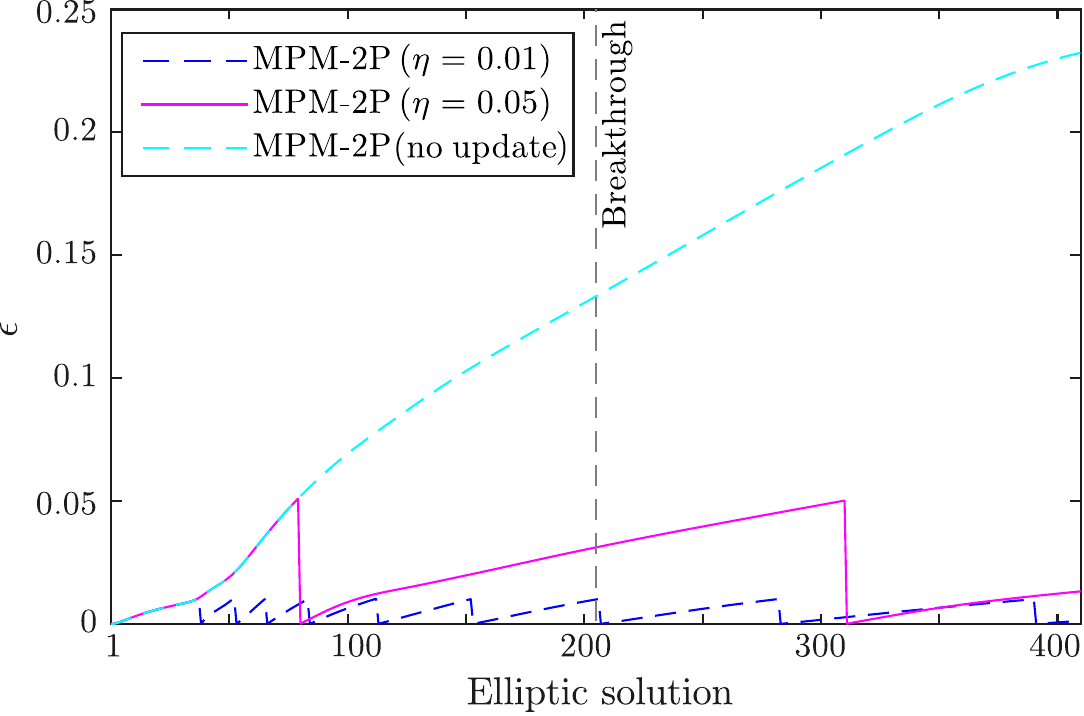}
	  \caption{Variation of $\epsilon$ throughout the simulation controlled by three different tolerance criteria: no BFs updates, $\eta=0.01$, and $\eta=0.05$.}
	\label{fig:MPM_5}	  
\end{figure}

In Fig. \ref{fig:MPM_4} we test the sensitivity of the tolerance $\eta$, by comparing the results yielded by MPM-2P for different choices of $\eta$ with the results obtained by the $a$MRCM-PBS. By changing from $\eta = 0.01$ to $\eta=0.05$, one can see a  decrease in the number of updates of the BFs in the same time frame from about 10 to 3. Note the error is still well bounded around $10^{-1}$ even with such low number of updates. We also turned off BFs updates, that shows an increase in the error, that is still well behaved, tending to converge to a value around $5\times 10^{-1}$ for this problem. 
The variation of $\epsilon$ throughout the simulation can be found in Fig. \ref{fig:MPM_5}, where we point out the tolerance criterion controlling its values.


\begin{figure}[htbp]
    \centering
	\includegraphics[scale=0.61]{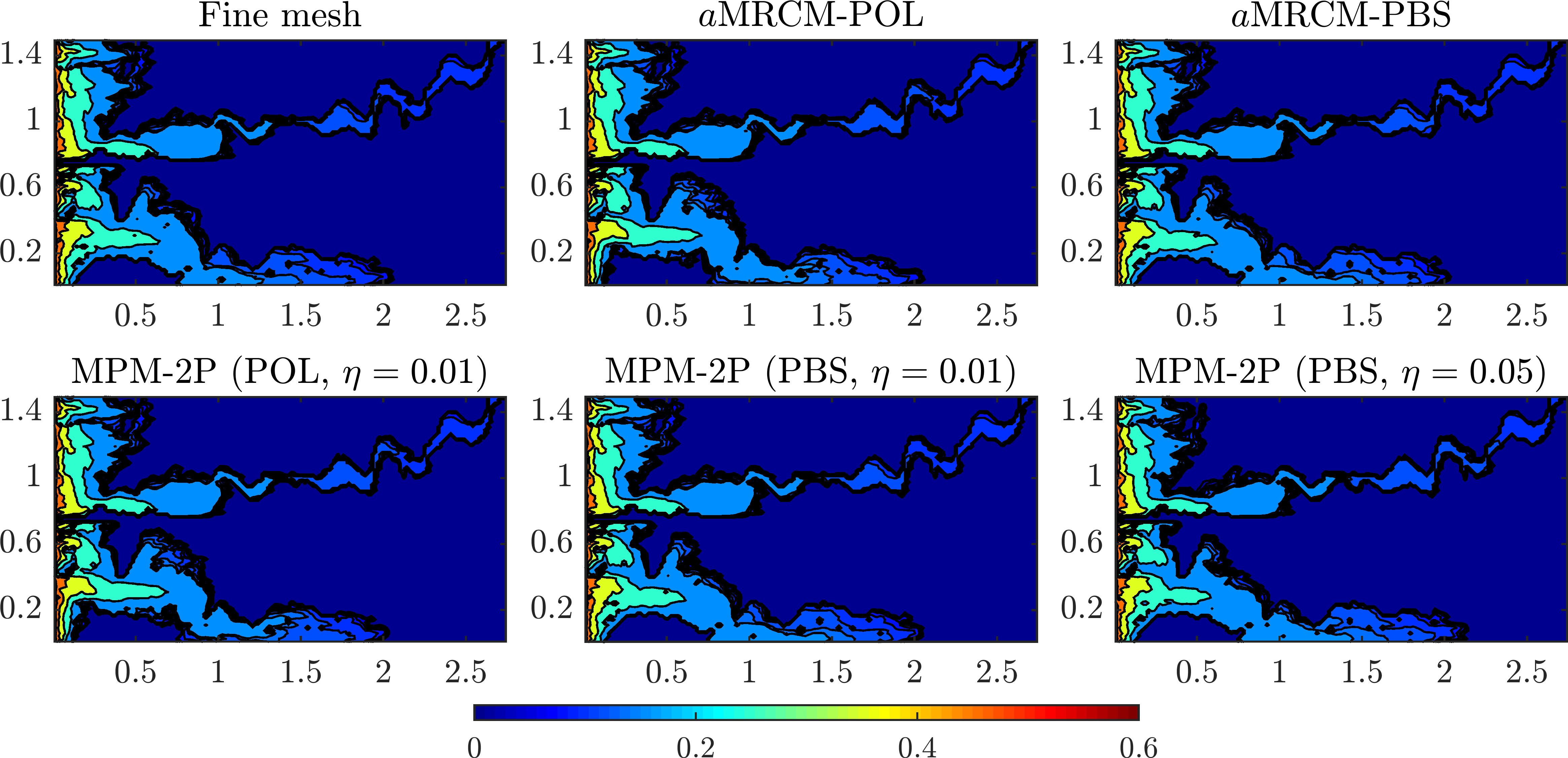}
	    \caption{Saturation profiles at breakthrough time $T_{\text{PVI}}=0.06$ (that corresponds to $206$ elliptic solutions), comparing the methods  $a$MRCM-POL, $a$MRCM-PBS, MPM-2P-POL ($\eta=0.01$), MPM-2P-PBS ($\eta=0.01$) and MPM-2P-PBS ($\eta=0.05$).}
	 \label{fig:MPM_6}
\end{figure}

To close this example, we compare the methods in terms of their saturation profiles in Fig. \ref{fig:MPM_6}. 
We show saturation plots for the fine mesh, $a$MRCM-POL, $a$MRCM-PBS, and the MPM-2P considering $\eta=0.01$ (combined with the $a$MRCM-POL and $a$MRCM-PBS) and $\eta=0.05$ (combined with the $a$MRCM-PBS). The profiles are taken at time $T_{\text{PVI}}=0.06$, that corresponds to the breakthrough time, i.e. $206$ elliptic solutions. 
Note that the improvement provided by the choice of the interface spaces base on physics (-PBS) over the $a$MRCM carry on to the MPM-2P as well. This accuracy is maintained even when the tolerance is relaxed to $\eta = 0.05$, which further decreases the number of BFs updates.

The Relative Cost Reduction attained by the MPM-2P for this numerical experiment can be found on Table \ref{table:Rocha_field}. These remarkable results show how much we can save by not recomputing all BFs every elliptic time step when solving two-phase flows through such heterogeneous media, and still keeping the solutions as accurate as those obtained by direct use of sophisticated multiscale mixed methods, with unprecedented reduction of the computational cost.


\begin{table}
\centering
\caption{Relative Cost Reduction attained by the MPM-2P using $a$MRCM-POL and $a$MRCM-PBS with different values of the tolerance $\eta$.} 
\begin{tabular}{c c}
\hline 
Method and tolerance & RCR \\ 
\hline 
$a$MRCM-POL $\eta=0.01$ &  84.06\% \\ 
$a$MRCM-PBS $\eta=0.01$ & 84.99\% \\ 
$a$MRCM-PBS $\eta=0.05$ & 87.13\% \\ 
$a$MRCM-PBS no updates & 88.43\% \\ 
\hline 
\end{tabular} 
\label{table:Rocha_field}
\end{table}

\subsection{A fractured permeability field}\label{subsec:A fractured permeability field}

\begin{figure}[htbp]
\centering
\includegraphics[scale=0.6]{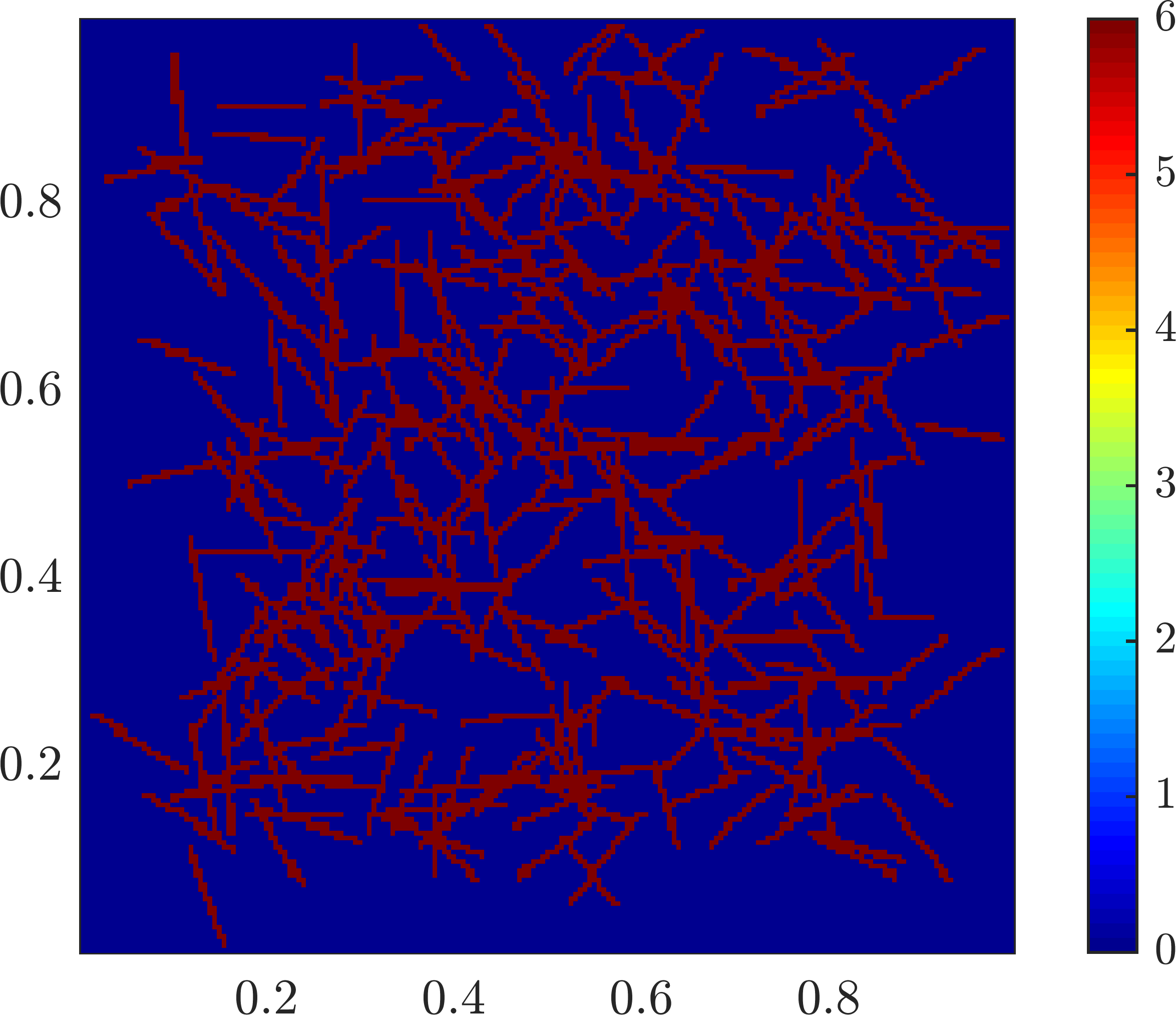}
\caption{Fractured permeability field (log-scaled).}
\label{fig:fractured}
\end{figure}

Another challenging test for the  MPM-2P is the fractured permeability field illustrated in Fig. \ref{fig:fractured}. The domain is set to $\Omega=[0,1]\times[0,1]$, with $200\times 200$ fine grid cells and a domain decomposition of $10\times10$ subdomains. The flow setup is about the same as in the previous experiment. The $a$MRCM-PBS is used to deal with the high-permeable fractures.

\begin{figure}[htbp]
\centering
\includegraphics[scale=0.44]{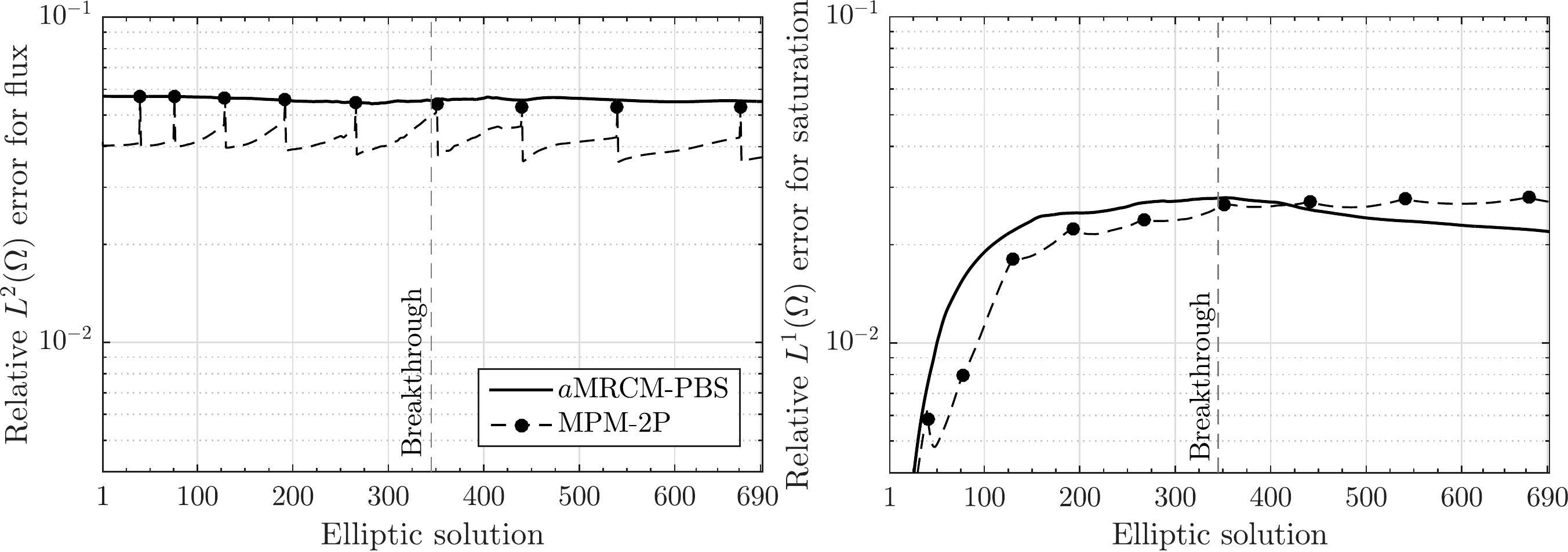}
\caption{Relative $L^2(\Omega)$ errors for flux (left) and relative $L^1(\Omega)$ errors for saturation (right). We consider physics-based interface spaces and compare the $a$MRCM-PBS and the MPM-2P (with $\eta=10^{-2}$). The nodes on each of the MPM-2P curves indicate the times when the BFs are updated. Note that the precision of the $a$MRCM-PBS and the MPM-2P are comparable.}
\label{fig:fractured_errors}
\end{figure}

Figure \ref{fig:fractured_errors} shows the relative $L^2(\Omega)$ errors for flux (left) and relative $L^1(\Omega)$ errors for saturation (right) computed with respect to the reference fine grid solution. The breakthrough time for the fine grid solution is indicated by a vertical dashed line.
We consider the tolerance of $\eta=10^{-2}$ for the updates of the BFs in the MPM-2P algorithm (represented by the nodes). Note that the set of BFs was computed 10 times by the MPM-2P (the initial set plus 9 updates), that is significantly less than the total of $692$ full updates required by the $a$MRCM-PBS. This reduction represents a significant reduction for the computational cost, while the accuracies of the approximations are comparable.

\begin{figure}[htbp]
\centering
\includegraphics[scale=0.69]{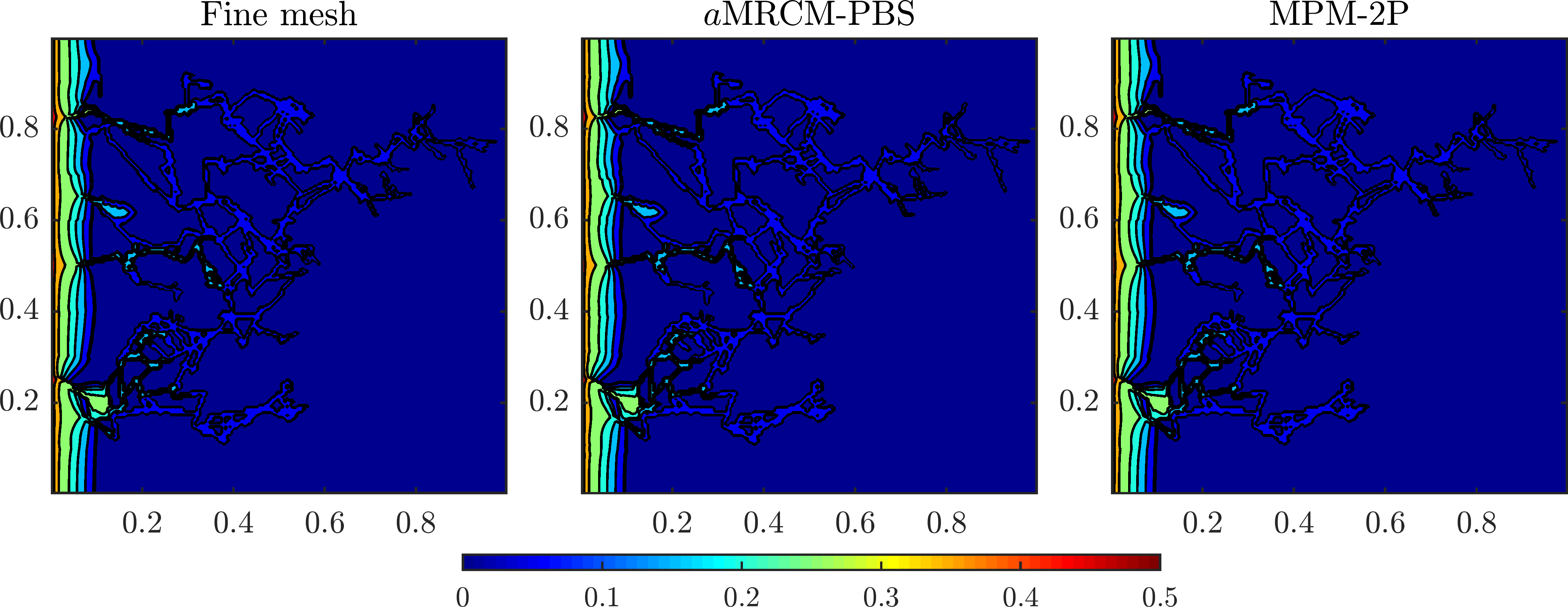}
\caption{Saturation profile at breakthrough time $T_{\text{PVI}}=0.03$, that corresponds to the elliptic time step number $346$. Left to right: fine grid solution; $a$MRCM-PBS; MPM-2P using the $a$MRCM-PBS for the updates with $\eta=0.01$.}
\label{fig:fractured_sat}
\end{figure}

A comparison of the saturation profiles at the breakthrough time $T_{\text{PVI}}=0.03$, that corresponds to $346$ elliptic solutions, is shown in Fig. \ref{fig:fractured_sat}. 
We show saturation maps for the fine mesh, $a$MRCM-PBS, and the MPM-2P combined with the $a$MRCM-PBS. 
Note that both approximations are closely related to the reference solution, being the MPM-2P approximation significantly less expensive than the $a$MRCM-PBS one. 
The Relative Cost Reduction attained by the MPM-2P for this  experiment is $RCR=88.45\%$.

\subsection{Water-oil finger growth in a homogeneous medium}\label{subsec:Water-oil finger growth in homogeneous media}

\begin{figure}[htbp]
    	\centering
	\includegraphics[scale=0.73]{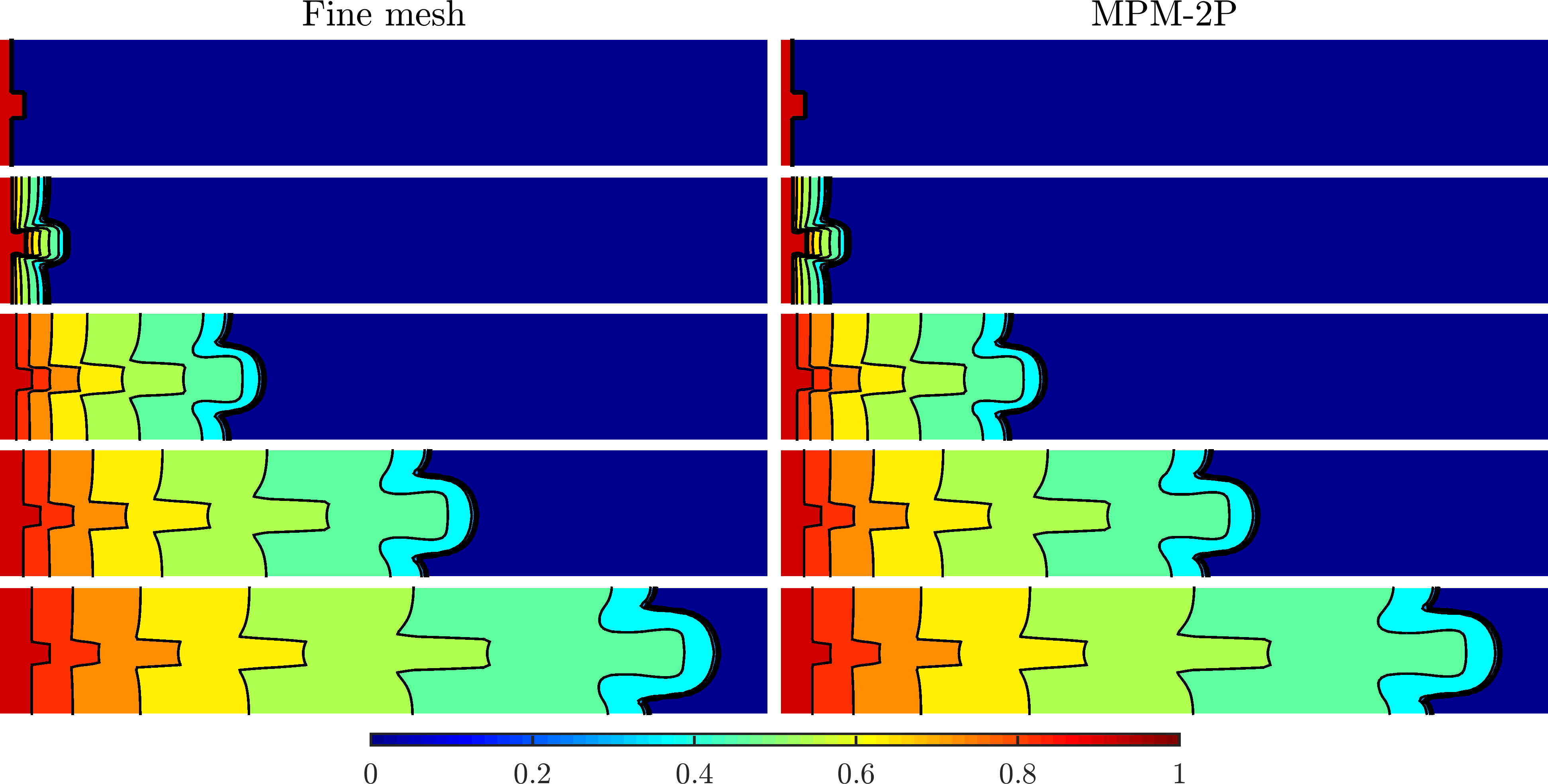}
	    \caption{Saturation evolution for the Riemann problem with a small perturbation of the initial water-oil interface at the center of the domain. We show the fine grid solution (left) and the MPM-2P approximation (right) at times $T_{\text{PVI}}=0.00$,\ $0.03$,\ $0.19$,\ $0.39$,\ $0.66$ (corresponding to $1$,\ $100$,\ $600$,\ $1100$, and $1600$ elliptic solutions, respectively), from top to bottom.}
	\label{fig:MPM_7}	    
\end{figure}

In this last experiment, we test the MPM-2P with an unstable interface of oil-water in a homogeneous medium. We consider an injection of water at the left boundary of the domain $\Omega=[0,3]\times[0,1/2]$, with an initial front fully saturated of water at the left and filled with oil at the right. The water front has a small perturbation at the center of the channel, as shown in Fig. \ref{fig:MPM_7}. This is a 2D Riemann problem with physical instabilities, similar to the studied in \cite{glimm1981numerical}, where the authors have shown that a finger grows at the center of the channel and evolves in time. Here we consider the same viscosity ratio $M=4$, that is close to the critical value for unstable flows ($M \approx 2.657$) as described in \cite{furtado2003crossover}. The boundary conditions considered are no-flow at the top and bottom along with an imposed pressure $p = 0$ on the left and $p = -10^{4}$ on the right boundaries. Furthermore, no source terms are considered.

Figure \ref{fig:MPM_7} shows the evolution of the saturation for times $T_{\text{PVI}}=0.00$, $0.03$,  $0.19$, $0.39$,  $0.66$ (corresponding to $1$, $100$, $600$, $1100$, and $1600$ elliptic solutions, respectively), from top to bottom. The reference fine grid solution (left) and the MPM-2P approximation (right) are compared in this figure. 
The MPM-2P uses a domain decomposition of $15\times 5$ subdomains, with $20\times10$ cells into each one.
For the update of BFs we use the MRCM with $\bar H = h$ and constant interface spaces, therefore, the solution obtained by the MRCM with these parameters is the same as the undecomposed case (see \cite{guiraldello2018multiscale}), which means that we do not have inaccuracies coming from the domain decomposition by the MRCM in this comparison. Moreover, inaccuracies associated with the operator splitting are also reduced by setting $C_n=1,\ \forall n$.
The updates of the BFs in the MPM-2P are performed according to the chosen tolerance of $\eta=10^{-2}$. 

\begin{figure}[htbp]
    	\centering
	\includegraphics[scale=0.65]{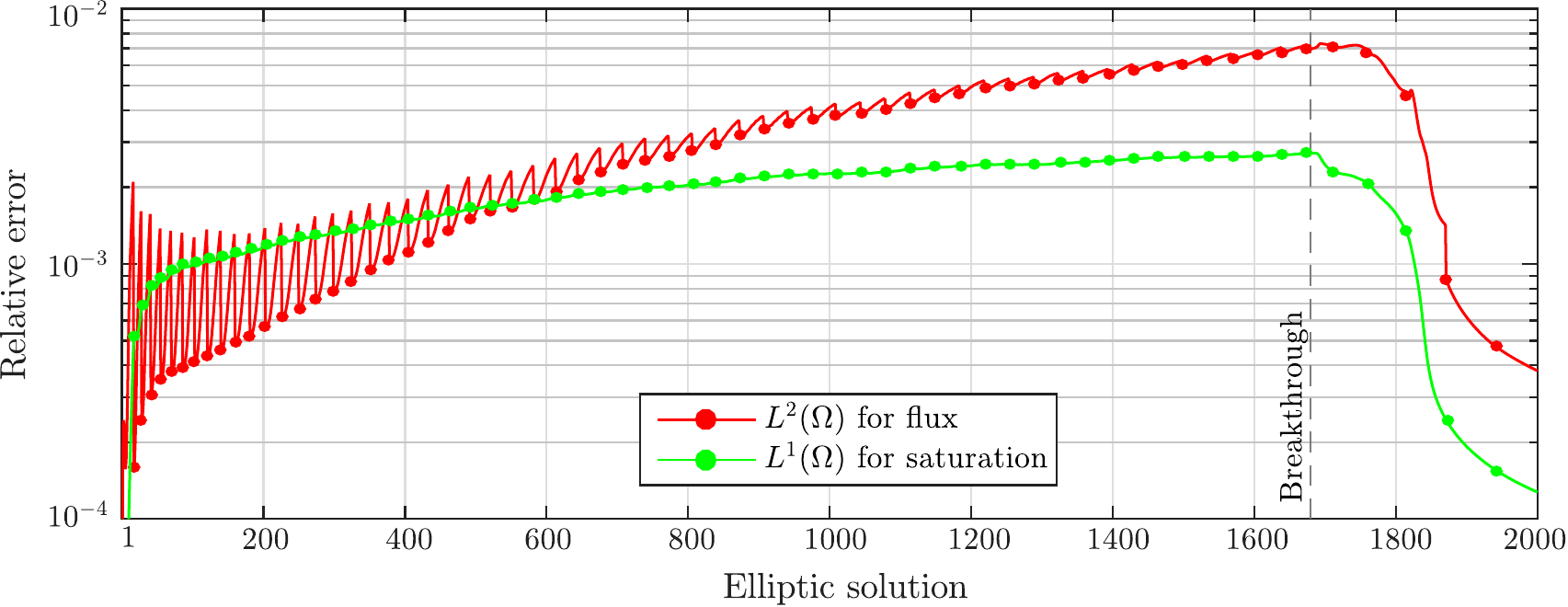}
	    \caption{Relative $L^2(\Omega)$ error for flux and relative $L^1(\Omega)$ error for saturation, obtained by the MPM-2P for the Riemann problem with a small perturbation at the center of the domain. The nodes indicate the times when BFs were updated.}
	\label{fig:MPM_8}
\end{figure}

The relative $L^2(\Omega)$ error for flux and relative $L^1(\Omega)$ error for saturation as functions of the number of elliptic solutions can be found in Fig. \ref{fig:MPM_8}. A total of $T_e=2000$ elliptic solutions were performed, and the MPM-2P solution required only $64$ updates of the BFs (indicated by the nodes in that figure). The trend of quickly increasing errors until the breakthrough time (illustrated by a dashed line) was controlled by the updates of the BFs. This illustrates how challenging is this problem and how it is well handled by the MPM-2P, with a rapid drop on flux error after every update. Even for this complex problem with physical instabilities, the MPM-2P decreases the number of full updates of the set of BFs from $2000$ to $64$, yielding a Relative Cost Reduction of  $RCR=94.62\%$. These results confirm the great potential the MPM-2P to reduce drastically the computational cost of two-phase flow simulations, without loss of accuracy, being suitable for any physically-challenging incompressible two-phase subsurface flow problem.

\section{Conclusion and future work}

In this work, we introduced and tested the Multiscale Perturbation Method for two-phase flows in porous media (MPM-2P). 
We formulate a modified operator splitting method, where we replace full updates of local solutions by reusing basis functions computed by a multiscale mixed method 
(such as the Multiscale Robin Coupled Method - MRCM) at an earlier time of the simulation. The reuse of multiscale basis functions is guaranteed by using perturbation 
theory to write suitable local problems, drastically reducing the computational cost of  multiscale mixed methods.

Our numerical results show an exceptional reduction in the computational cost of the simulation of two-phase flows in challenging permeability fields. The MPM-2P 
can improve significantly the efficiency 
of an operator splitting method for two-phase flows, without loss of accuracy. The numerical examples show that water breakthrough can be simulated with very few 
updates of the MRCM set of basis functions. 
The errors produced by the MPM-2P are comparable, and in most cases smaller, to the typical values of error attained by multiscale mixed methods.  We remark that any 
multiscale mixed method can be used for the updates of the basis functions of the MPM-2P formulation in a straightforward manner.

The implementation of the new method in multi-core and multi-GPU devices and its application to the sequential implicit solution of multiphase flows are currently being considered by the authors and their collaborators. Moreover, the use of MPM-2P in accelerating Markov chain Monte Carlo methods for uncertainty quantification of subsurface flows is a promising research topic and is also being investigated by the authors.


\section*{Acknowledgements}

F. F. Rocha, F. S. Sousa and F. Pereira acknowledge the financial support received from Brazilian oil company Petrobras grant 2015/00400-4, and from the S\~ao Paulo Research Foundation FAPESP, CEPID-CeMEAI grant 2013/07375-0; This study was also funded in part by Brazilian government agencies CAPES (Finance Code 001) and CNPq; F. S. Sousa was funded in part by CNPq grant 310990/2019-0. 

\bibliography{MPMRef1,MPMRef2}
\bibliographystyle{elsarticle-num} 
\end{document}